\documentclass[11pt,letterpaper]{article}
\usepackage[utf8]{inputenc}

\usepackage{pdflscape}
\usepackage{appendix} 
\usepackage{amsmath,bm}
\usepackage{amsfonts}
\usepackage{graphicx}
\usepackage{ragged2e}
\usepackage{graphicx}
\usepackage{multicol}
\usepackage{placeins}
\usepackage[hidelinks]{hyperref}
\usepackage{xr-hyper}
\usepackage{xcolor}
\usepackage{tabularray}
\UseTblrLibrary{siunitx}
\usepackage{graphicx}
\usepackage{amsmath,bm}
\usepackage{amsmath}
\usepackage{graphicx}
\usepackage{amsfonts}
\usepackage{xcolor}

\usepackage[table]{xcolor}
\usepackage[left=2.5cm,right=2.5cm,top=1.5cm,bottom=1.5cm]{geometry}

\usepackage{tikz}
\usepackage{authblk}
\usepackage{pgfplots}
\pgfplotsset{compat=1.17}

\usepackage{xr} % or \usepackage{xr-hyper} if using hyperref

\begin{document}
\title{Representation of tensor functions using lower-order structural tensor set: three-dimensional theory}
\author{Mohammad Madadi and Pu Zhang\thanks{Corresponding author, e-mail: pzhang@binghamton.edu}}
\affil{\textit{Department of Mechanical Engineering, State University of New York at Binghamton, Binghamton, NY 13902, USA}}
\date{}
\maketitle

\begin{justify}
\begin{abstract}
The representation theory of tensor functions is a powerful mathematical tool for constitutive modeling of anisotropic materials. A major limitation of the traditional theory is that many point groups require fourth- or sixth-order structural tensors, which significantly impedes practical engineering applications. Recent advances have introduced a reformulated representation theory that enables the modeling of anisotropic materials using only lower-order structural tensors (i.e., second-order or lower). Building upon the reformulated theory, this work establishes the representations of tensor functions for three-dimensional centrosymmetric point groups. For each point group, we propose a lower-order structural tensor set and derive the representations of tensor functions explicitly. For scalar-valued and second-order symmetric tensor-valued functions, our theory is indeed applicable to all three-dimensional point groups because their representations are determined by the corresponding centrosymmetric groups. The representation theory presented here is broadly applicable for constitutive modeling of anisotropic materials.
\end{abstract}

\section{Introduction}\label{sec:introduction}
Constitutive modeling is central to continuum mechanics and materials modeling in engineering and materials science \cite{ottosen2005mechanics}. Constitutive models capture a wide range of mechanical and multiphysical behaviors of materials including stress–strain relationships, yield surfaces, failure criteria, thermal and electrical properties, and mechano-physical behaviors. Constitutive laws are usually described using scalar- or tensor-valued functions with multiple arguments including field and state variables. In addition, material symmetry or anisotropy is incorporated in constitutive laws through structural tensors (i.e., anisotropic tensors) designated for each point group \cite{zheng_theory_1994}.  

The representation theory of tensor functions \cite{zheng_theory_1994, Boehler_applications_1987, apel2004approaches} is a powerful mathematical tool that provides general forms of constitutive laws consistent with frame-indifference and material symmetry principles. It was established in the mid-20th century by Rivlin  \cite{rivlin_large_1948, rivlin_further_1955}, Pipkin \cite{pipkin_formulation_1959}, and Noll \cite{noll_representations_1970}. These pioneers developed representations of isotropic tensor functions, primarily for isotropic materials. Thereafter, Wang \cite{wang_representations_1969, wangb_representations_1969}, Smith \cite{smith_isotropic_1971}, and Boehler \cite{boehler_irreducible_1977} derived isotropic scalar-, vector-, and tensor-valued functions of vectors and 2nd-order tensors. Later on, the representation theory was generalized to anisotropic tensor functions by Boehler and Liu \cite{boehler_simple_1979, liu_representations_1982}, and then Spencer and Betten \cite{spencer_formulation_1982, betten_irreducible_1987}. Their approach was to transform an anisotropic tensor function into an extended isotropic one by including structural tensor arguments. These structural tensors characterize material symmetry and are invariant under any symmetry operation of the point group \cite{zheng_theory_1994}.  The structural tensors for all two-dimensional (2D) and three-dimensional (3D) point groups are reviewed by Zheng \cite{zheng_theory_1994}. Despite its theoretical elegance, the practical application of this anisotropic representation theory is rather limited. Firstly, the theory provides only the general mathematical forms, requiring researchers to determine specific functions either empirically or through trial-and-error. Secondly, a major obstacle is that most point groups involve higher-order (i.e., 3rd-order or higher) structural tensors, which complicate and often preclude practical modeling. In practice, constitutive modeling with higher-order structural tensors is rarely feasible. Most crystalline point groups require higher-order structural tensors and their constitutive modeling remains largely unexplored. To circumvent the obstacle of higher-order structural tensors, Xiao \cite{xiao_isotropic_1996} introduced lower-order structural tensor functions (i.e., non-constant structural tensors) to derive the representations of anisotropic tensor functions. Key results of the representations are reported in following works \cite{xiao1996minimal, xiao1998scalar}.

Alternatively, to circumvent the obstacle of higher-order structural tensors, Man and Goddard reformulated the representation theory in 2018 \cite{man2018remarks}, enabling the exclusive use of lower-order structural tensors. Unlike the original theory of Boehler and Liu \cite{boehler_simple_1979, liu_representations_1982}, which requires structural tensors to be invariant under all symmetry operations of a point group, the Man-Goddard reformulation relaxes this requirement while imposing additional symmetry constraints afterwards. This reformulation enables the constitutive modeling of anisotropic materials using only lower-order structural tensors. In their work \cite{man2018remarks}, Man and Goddard provided illustrative examples demonstrating the reformulation, but did not fully establish the representations for all point groups. Based on the Man-Goddard reformulation, our recent work \cite{madadi2025representation} introduced a new concept "structural tensor set" and established the representations of tensor functions for all 2D point groups using lower-order structural tensor sets. Besides, our recent work also presented comprehensive review and discussions on the original and reformulated representation theories. Despite the progress for 2D point groups, a critical knowledge gap remained regarding the reformulated representation theory for 3D point groups. Very recently, Man and Goddard \cite{man2026further} provided a rigorous mathematical framework for this reformulated representation theory. They argued that a point group may be characterized by multiple structural tensor sets, with the results of Boehler-Liu \cite{boehler_simple_1979, liu_representations_1982} and Madadi \cite{madadi2025representation} as special cases. The present work aims to establish detailed representations of tensor functions for all 3D centrosymmetric point groups, including 11 Laue groups and 3 continuous ones. Our previous work discovered that, for a given point group, the representations of scalar-valued and 2nd-order symmetric tensor-valued functions are determined by its corresponding centrosymmetric group (e.g., Laue group) \cite{madadi2025representation}. Hence, we limit our theory to centrosymmetric point groups in this work. As long as scalar-valued and 2nd-order symmetric tensor-valued functions are of interest, the presented theory is applicable to all 3D point groups (i.e., 32 crystalline point groups and 7 continuous ones) because one only needs to find their corresponding centrosymmetric groups and the associated representations in this work.

The representation of anisotropic scalar- and tensor-valued functions has broad applications in engineering and materials science \cite{ottosen2005mechanics, zheng_theory_1994,Boehler_applications_1987}. For scalar-valued functions, the theory developed here can be used to model hyperelastic strain energy functions of elastomers, soft composites \cite{ma2025precisely, Rahman2025}, and biological tissues \cite{chagnon2015hyperelastic}, as well as yield and failure criteria for materials \cite{ottosen2005mechanics,Boehler_applications_1987}. For tensor-valued functions, the representation theory provides a basis for modeling mechanical, physical, and mechano-physical properties \cite{madadi_finite_2024}, including stress–strain relations, dielectric properties, and conductivity tensors \cite{nguyen_modeling_2024}. It should be emphasized that the present work establishes only the general forms of such scalar- and tensor-valued functions. Specific constitutive laws must still be constructed and fitted to experimental or simulation data. 

\sloppy The remainder of this paper is organized as follows. In Sections \ref{sec:2}, we revisit the preliminaries of representation theory for tensor functions. Section \ref{sec:3} introduces the proposed lower-order structural tensor sets for all 3D centrosymmetric point groups. The complete representations of scalar- and 2nd-order symmetric tensor-valued functions are reported in Sections \ref{sec:C_i}–\ref{sec:Continuous}. We have tried our best to simplify the representations and remove redundant terms, although irreducibility is only desirable but not mandatory. Notably, six groups ($\mathcal C_i, \mathcal C_{2h}, \mathcal D_{2h}, \mathcal C_{\infty h}, \mathcal D_{\infty h}, \mathcal K_{h}$) possess only lower-order structural tensors and can therefore be treated using the original Boehler–Liu formulation. In contrast, eight groups ($\mathcal C_{4h}, \mathcal D_{4h}, \mathcal C_{3i}, \mathcal D_{3d}, \mathcal C_{6h}, \mathcal D_{6h}, \mathcal T_{h}, \mathcal O_{h}$) possess higher-order structural tensors and should employ the Man-Goddard reformulation together with our proposed lower-order structural tensor sets. Throughout this work, we adopt the Schoenflies notation for point groups; for other notation systems, readers may refer to \cite{de2012structure}. Finally, the appendix lists functional bases, tensor generators, symmetry operations, and useful matrices for reference.

\section{Preliminaries of representation theory} \label{sec:2}
In our previous work \cite{madadi2025representation}, we have provided a brief introduction to the representation theory of tensor functions in 2D space. The corresponding theory for 3D is similar. The major formulas are summarized below.  

Firstly, we consider the representation of isotropic scalar- and tensor-valued functions. Isotropic tensor functions are useful for the constitutive modeling of isotropic materials. In general, a scalar-valued isotropic tensor function $\psi(\mathbf v, \mathbf{A}, \mathbf{W})$ can be expressed as a function of the invariants $I_k$ \cite{ajm_theory_1971} of its arguments $\mathbf v$, $\textbf{A}$, and $\mathbf{W}$, as
\begin{align} \label{eq_iso_scalar_3D}
	\psi(\mathbf v, \mathbf{A}, \mathbf{W})=\psi(I_k)
\end{align}
where $\mathbf{v}$, $\mathbf{A}$, and $\mathbf{W}$ are sets of vectors, 2nd-order symmetric and skew-symmetric tensors, respectively. Herein, the complete set of invariants $I_k(k=1,2,...,r)$ is called functional basis (or integrity basis). The representation of an isotropic 2nd-order symmetric tensor-valued function $\mathbf{T}(\mathbf{v}, \mathbf{A}, \mathbf{W})$ is expressed as a linear combination of tensor generators $\mathbf G_i$, as  
\begin{align} \label{eq_iso_second_3D}
	\begin{array} {ll}
		\mathbf{T}(\mathbf{v}, \mathbf{A}, \mathbf{W})=\sum\alpha_i \mathbf G_i
	\end{array}
\end{align}
where $\alpha_i = \alpha_i(I_1,I_2,...,I_r)$ are scalar coefficient functions of the invariants $I_k$ governed by (\ref{eq_iso_scalar_3D}). Once the arguments $\mathbf v$, $\textbf{A}$, and $\mathbf{W}$ are provided, one can find the tensor generators $\mathbf{G}_i$ following methods and formulae presented in \cite{Boehler_applications_1987, zheng_theory_1994}. For an isotropic 2nd-order symmetric tensor-valued function $\mathbf{T}(\mathbf{v}, \mathbf{A}, \mathbf{W})$, the functional bases $I_k$ and tensor generators $\mathbf{G}_i$ can be obtained from Table \ref{Table functional basis_3D} and Table \ref{Table generators_3D} in the Appendix.

Secondly, we consider the representation of anisotropic scalar- and tensor-valued functions, which are needed for constitutive modeling of anisotropic materials. By introducing structural tensors $\mathfrak{M}$, Boehler \cite{boehler_simple_1979} and Liu \cite{liu_representations_1982} introduced isotropic tensor functions $\hat{\psi}(\mathbf v, \mathbf{A},\mathbf{W}, \mathfrak{M})$ and $\hat{\mathbf T} (\mathbf v, \mathbf{A},\mathbf{W},\mathfrak{M})$ for anisotropic materials, which are actually isotropic extension of anisotropic functions. The requirement is that the structural tensors $\mathfrak{M}$ must be invariant under any symmetry operation $\mathbf Q$ in the point group $\mathcal{G}$. General forms of the extended isotropic functions $\hat \psi$ and $\hat{\mathbf{T}}$ can be derived readily using (\ref{eq_iso_scalar_3D}) and (\ref{eq_iso_second_3D}), respectively. The challenge is that many point groups involve higher-order structural tensors, which hinder the wide applications of the Boehler-Liu formulation.  

In order to overcome the challenge of higher-order structural tensors, Man and Goddard \cite{man2018remarks} reformulated the representation theory of anisotropic tensor functions. In their reformulation, only lower-order structural tensors $\mathfrak{M}$ are needed. The representations are similar to Boehler and Liu above. However, an additional symmetry constraint should be imposed to the tensor functions \cite{man2018remarks, madadi2025representation}, as
\begin{align} \label{eq_constraint}
	\begin{array} {ll}
		{\hat \psi}(\mathbf v, \mathbf{A}, \mathbf{W}, \mathfrak M) = {\hat \psi}(\mathbf{v}, \mathbf{A}, \mathbf{W}, \langle\mathbf Q \rangle \mathfrak{M}) \\	
		\mathbf{\hat T}(\mathbf v, \mathbf{A}, \mathbf{W},\mathfrak M) = \mathbf{\hat T}(\mathbf{v}, \mathbf{A}, \mathbf{W}, \langle\mathbf Q \rangle\mathfrak{M}) 
	\end{array} ;\forall \mathbf Q \in \mathcal{G}^*  
\end{align}
where $\mathcal{G}^*$ denotes the group generators of the point group $\mathcal G$, and we \cite{madadi2025representation} have proven that only the group generators need to be considered in (\ref{eq_constraint}) rather than the whole point group. The orthogonal transformation operator $\langle \mathbf Q\rangle$ is defined after Zheng \cite{zheng_theory_1994}, as
\begin{align} \label{eq 2.7}
	\begin{array} {ll}
		\langle\mathbf Q \rangle \mathbf v = \mathbf Q \mathbf v=  \mathbf Q_{ij} \mathbf v_j \\
		\langle\mathbf Q \rangle \mathbf A = \mathbf{Q} \mathbf{A} \mathbf{Q}^T= \mathbf Q_{ip} \mathbf Q_{jq} \mathbf A_{pq}\\ 
		\langle\mathbf Q \rangle \mathbb A = \mathbf Q_{ip} \mathbf Q_{jq} \mathbf Q_{kr}...\mathbf Q_{mt} \mathbb A_{pqr...t}
	\end{array}
\end{align}
where $\mathbf v$, $\mathbf A $, and $\mathbb A$ are first-, second-, and higher-order tensors, respectively. 

Using the Man-Goddard reformulation, we have established the representation theory of tensor functions for all 2D point groups in a previous work \cite{madadi2025representation}. The purpose of the present work is to establish the representation theory of tensor functions for all 3D centrosymmetric groups (Figure \ref{fig_3D point groups latiice}), including 11 Laue groups and 3 continuous groups, which are the most useful groups for constitutive modeling of anisotropic materials.

\begin{figure}[h!]
	\centering\includegraphics[width=5in]{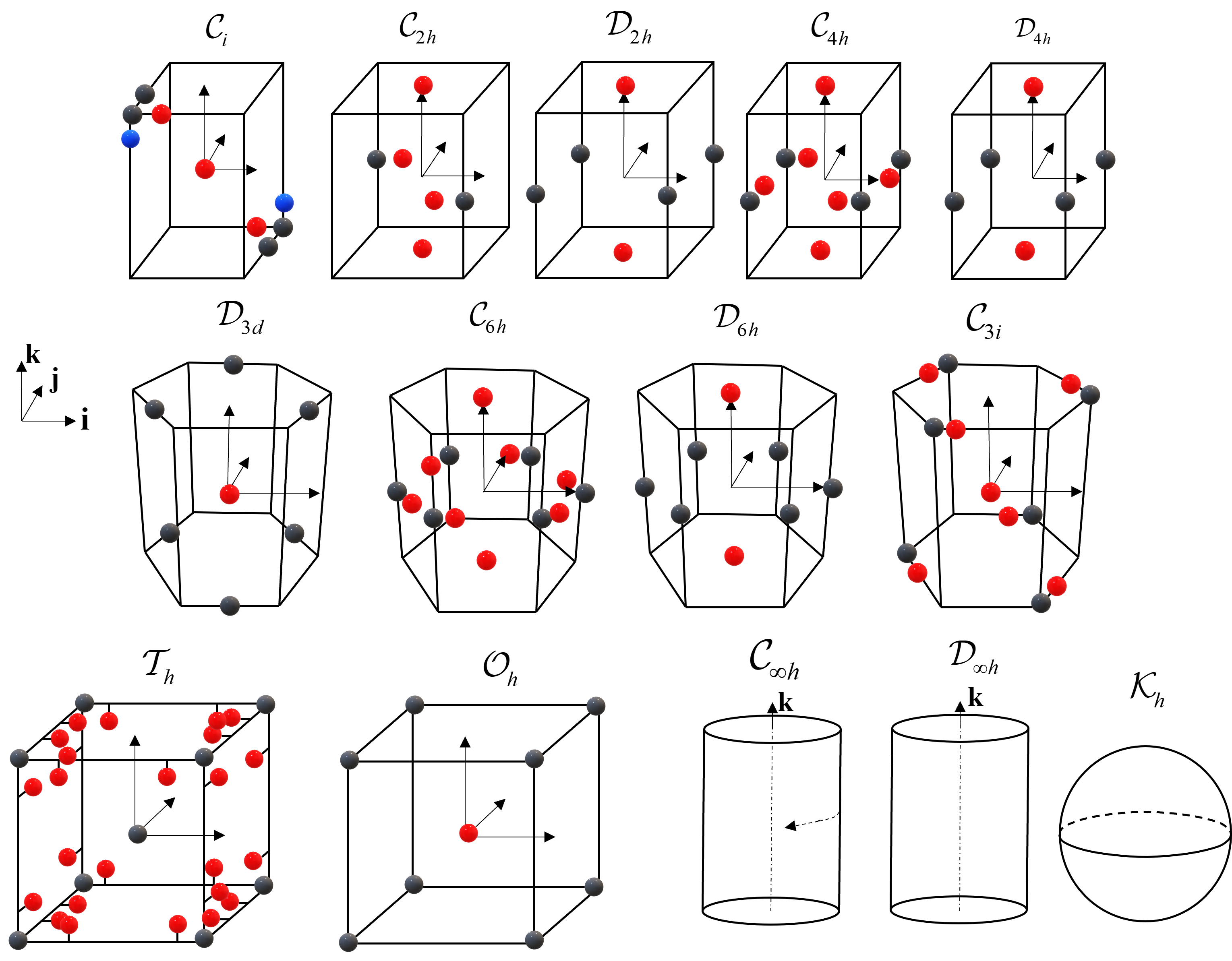}
	\caption{Graphical illustration of 3D centrosymmetric point groups: 11 Laue groups and 3 continuous groups.}
	\label{fig_3D point groups latiice}
\end{figure}
\vspace*{-5pt}

\section{Lower-order structural tensor set for 3D point groups} \label{sec:3}
In our previous work \cite{madadi2025representation}, we proposed the concept "structural tensor set" and provided specific lower-order structural tensors for each 2D point group. Very recently, Man and Goddard \cite{man2026further} proposed that a point group can be characterized by multiple structural tensor sets. In this work, we adopt only a single structural tensor set for simplicity purposes. Specifically, the structural tensor set $\{ \mathfrak{M}_i \}$ of a 3D point group $\mathcal{G}$ must satisfy two conditions below.
	\begin{align} \label{eq_structur_def_3D_1}
		{\mathcal G_s=\{ \mathbf Q \in \mathcal O (3)\ |\  \langle \mathbf Q\rangle \mathfrak M_i  = \mathfrak M_i  ;\ \ i=1,2,...,t\}
		}
	\end{align}
	\vspace{-12mm}
	\begin{align} \label{eq_structur_def_3D_2}
		{\{\langle \mathbf Q\rangle \mathfrak M_i \}= \{\mathfrak M_i \},\  \forall \mathbf Q \in \mathcal{G}
		}
	\end{align}
	In (\ref{eq_structur_def_3D_1}), $\mathcal G_s$ is a subgroup of $\mathcal G$ (i.e. $\mathcal{G}_s \le \mathcal{G}$) and $\mathcal{O}(3)$ denotes the 3D orthogonal group.  Ideally,  $\mathcal G_s$ is the largest subgroup of  $\mathcal G$ that can be characterized by lower-order structural tensors $ \mathfrak{M}_i $ following the traditional definition of Boehler and Liu \cite{zheng_theory_1994}. The idea here is that we can first find a set of lower-order structural tensors $\{ \mathfrak{M}_i \}$ to characterize some symmetry (i.e., $\mathcal{G}_s$) of the point group $\mathcal{G}$ using (\ref{eq_structur_def_3D_1}), then derive the representations using the Boehler-Liu formulation \cite{boehler_simple_1979, liu_representations_1982}, and finally impose additional symmetry constraints (\ref{eq_constraint}) for $\mathcal{G}$ following the Man-Goddard approach \cite{man2018remarks, madadi2025representation}. Moreover, in order to employ the Man-Goddard reformulation \cite{man2018remarks}, the whole structural tensor set $\{ \mathfrak{M}_i \}$ must be invariant under $\mathcal{G}$, which is satisfied by (\ref{eq_structur_def_3D_2}).  This new definition of structural tensor set in (\ref{eq_structur_def_3D_1})-(\ref{eq_structur_def_3D_2}) for 3D point groups can be generalized to 2D point groups in lieu of the original definition by Madadi \cite{madadi2025representation}. But the original definition for 2D point groups \cite{madadi2025representation} is too restrictive for 3D point groups in practice.

For each point group, the structural tensor set $\{\mathfrak{M}_i\}$ is non-unique. One has to devise a structural tensor set that is convenient to use and provides compact mathematical formulae. It usually takes laborious work and a lengthy trial-and-error process to find a complete structural tensor set. Generally, one can choose typical high symmetry directions, lines, and planes to construct the lower-order structural tensors. In what follows, we will take a 3D point group $\mathcal D_{4h}$ in Figure \ref{fig_3D various vectors} as an example and present three different approaches to construct its structural tensor set.

\begin{itemize}
	\item \textbf{Approach I}: Consider a vector $\mathbf v_1 = \frac{\sqrt 2}{2} \mathbf{i} + \frac{\sqrt 2}{2} \mathbf{j} + \mathbf{k} $ illustrated in Figure \ref{fig_3D various vectors} (a), which is within a high symmetry plane. We can first define a 2nd-order structural tensor $\mathbf M_1 = \mathbf v_1 \otimes \mathbf v_1$. By applying the point group generators provided in Table \ref{Table Laue generators}, the remaining structural tensors are found as $\mathbf{M}_2= \mathbf C_4 \mathbf M_1 \mathbf C_4^T$, $\mathbf{M}_3=\mathbf C_4 \mathbf M_2 \mathbf C_4^T$, and $\mathbf{M}_4=\mathbf C_{2x} \mathbf M_1 \mathbf C_{2x}^T$. 
	We can further verify that $\{\mathbf{M}_1, \mathbf{M}_2, \mathbf{M}_3, \mathbf{M}_4 \}$ form a complete structural tensor set. In this case, $\mathcal{G}_s=\mathcal{C}_i$ in (\ref{eq_structur_def_3D_1}).
	
	\item \textbf{Approach II}: Consider two orthonormal vectors $\mathbf v_1 ' = \frac{\sqrt 2}{2} {\mathbf i} + \frac{\sqrt 2}{2} {\mathbf j}$ and $\mathbf v_2 ' = {\mathbf k} $ illustrated in Figure \ref{fig_3D various vectors} (b). Herein, both vectors are along high symmetry axes. We can start with two structural tensor $\mathbf M_1 = \mathbf v_1 ' \otimes \mathbf v_1 '$ and $\mathbf M_2 = \mathbf v_2 ' \otimes \mathbf v_2 '$. By applying the group generators in Table \ref{Table Laue generators}, one extra structural tensor $\mathbf{M}_3 = \mathbf C_4 \mathbf M_1 \mathbf C_4^T$ is found. We can verify that $\{\mathbf{M}_1, \mathbf{M}_2, \mathbf{M}_3 \}$ form a complete structural tensor set. In this case, $\mathcal{G}_s=\mathcal{D}_{2h}$ in (\ref{eq_structur_def_3D_1}).
	
	\item \textbf{Approach III}: Consider three orthonormal vectors $\mathbf v_1 '' = \mathbf i, \mathbf v_2 '' = \mathbf j, \text{ and }\mathbf v_3 '' = \mathbf k$ illustrated in Figure \ref{fig_3D various vectors} (c). These three vectors are all along high symmetry axes and coincide with the coordinate axes. We can define three structural tensors as $\mathbf M_1 = \mathbf v_1 '' \otimes \mathbf v_1 '',\  \mathbf M_2 = \mathbf v_2 '' \otimes \mathbf v_2 '', \ \text{and }\mathbf M_3 = \mathbf v_3 '' \otimes \mathbf v_3 '' $. Further, we can verify that $\{\mathbf{M}_1, \mathbf{M}_2, \mathbf{M}_3 \}$ form a complete structural tensor set. In this case, $\mathcal{G}_s=\mathcal{D}_{2h}$ in (\ref{eq_structur_def_3D_1}).
	
\end{itemize}

\begin{figure} [h!]
	\centering\includegraphics[width=4in]{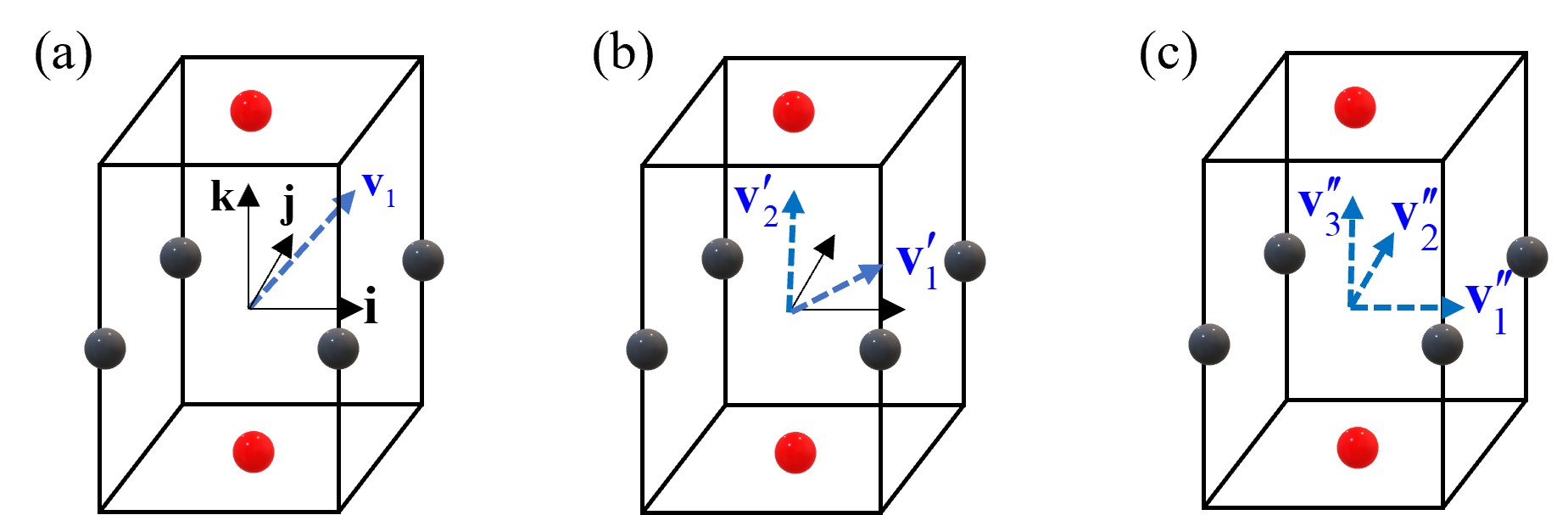}
	\caption{Illustration of various vectors used to define structural tensor set for point group $\mathcal D_{4h}$}.
	\label{fig_3D various vectors}
\end{figure}

It is generally preferable to select a structural tensor set $\{ \mathfrak M_i \} $ with few members and simple expressions because it would simplify the representations. For example, the Approach III above is strongly recommended for the point group $\mathcal{D}_{4h}$ for its simplicity.   

The structural tensors of 3D centrosymmetric point groups are presented in Table \ref{Table structural tensors_3D}. For the 6 groups ($\mathcal C_i, \mathcal C_{2h}, \mathcal D_{2h}, \mathcal C_{\infty h},\mathcal D_{\infty h}, \mathcal K_{\infty} $) with lower-order structural tensors, we simply adopt these tensors provided by Zheng \cite{zheng_theory_1994}. For these groups, the Boehler-Liu formulation should be used. In contrast, for the 8 point groups ($\mathcal C_{4h}, \ \mathcal D_{4h}, \ \mathcal C_{6h}, \ \mathcal D_{6h},\ \mathcal T_h, \ \mathcal O_h, \ \mathcal D_{3d}, \ \mathcal C_{3i}$) with higher-order structural tensors given by Zheng, we propose lower-order structural tensor sets for them. For these 8 groups, the Man-Goddard reformulation should be used. The representations of scalar- and 2nd-order symmetric tensor-valued functions for all 3D centrosymmetric point groups are presented in Sections \ref{sec:C_i} to \ref{sec:Continuous}. 

\begin{table}[!h]
	\caption{ Structural tensors for 3D centrosymmetric point groups}
	\label{Table structural tensors_3D}
	\centering
	\begin{tabular}{|l|c|c|}
		\hline
		\centering
		System (Point group) & Zheng's structural tensors \cite{zheng_theory_1994} & Proposed structural tensor set \\ [1.0ex] \hline
		Triclinic $(\mathcal C_i) $  & \footnotemark[1]$\bm \varepsilon \mathbf i$, $\bm \varepsilon \mathbf j$, $\bm \varepsilon \mathbf k$ & $\bm \varepsilon \mathbf i$, $\bm \varepsilon \mathbf j$, $\bm \varepsilon \mathbf k$ \\ [1ex]
		Monoclinic $\mathcal{(C}_{2h})$   & $\mathbf{P}_2, \bm{\varepsilon} \mathbf k$ & $\mathbf{P}_2, \bm{\varepsilon} \mathbf k$ \\ [1ex] 
		Orthorhombic $(\mathcal{D}_{2h})$  & $\mathbf{P}_2 $ & $\mathbf{P}_2 \ (\text{or } \mathbf M_1, \mathbf M_2, \mathbf M_3 ) $\\[1ex] 
		Tetragonal $(\mathcal{C}_{4h})$    &$\mathbb{P}_4, \bm{\varepsilon} \mathbf k $ &$\mathbf M_1, \mathbf M_2, \mathbf M_3, \bm\varepsilon \mathbf k $\\
		Tetragonal $(\mathcal{D}_{4h})$    &$\mathbb{P}_4 $ &$\mathbf M_1, \mathbf M_2, \mathbf M_3$\\
		[1ex] 
		Trigonal ($\mathcal{C}_{3i}  $) & $\mathbf k \otimes \mathbb{P}_3, \bm{\varepsilon} \mathbf k $ & $\mathbf T_1, \mathbf T_2, \mathbf T_3, \bm{\varepsilon} \mathbf k$  \\
		Trigonal ($\mathcal{D}_{3d}  $) & $\mathbf k \otimes \mathbb{P}_3$ & $\mathbf D_1, \mathbf D_2, \mathbf D_3$\\ [1ex]
		Hexagonal ($ \mathcal{C}_{6h}   $) & $\mathbb{P}_6, \bm{\varepsilon} \mathbf k$ & $\mathbf H_1,\mathbf H_2,\mathbf H_3, \bm{\varepsilon}\mathbf k $  \\
		Hexagonal ($ \mathcal{D}_{6h}   $) & $\mathbb{P}_6$ & $\mathbf H_1,\mathbf H_2,\mathbf H_3$  \\ [1ex]
		Cubic ($ \mathcal{T}_{h}$)  &  $\mathbb T_h$ &  $\mathbf M_1, \mathbf M_2, \mathbf M_3$ \\
		Cubic ($ \mathcal{O}_{h}$)   & $\mathbb O_h$ & $\mathbf M_1, \mathbf M_2, \mathbf M_3$ \\ [1ex]
		Cylindrical ($ \mathcal{C}_{\infty h}$)   & $\bm\varepsilon \mathbf k$ & $\bm\varepsilon \mathbf k$ \\
		Cylindrical ($ \mathcal{D}_{\infty h}$)   & $\mathbf k \otimes \mathbf k$ & $\mathbf k \otimes \mathbf k$ \\ [1ex]
		Spherical ($ \mathcal{K}_{h}$)   & $\mathbf I$ & $\mathbf I$ \\ [1.2ex]
		\hline
	\end{tabular}
	\footnotetext[1]{$\bm{\varepsilon}$ is the third-order permutation tensor in 3D space.}
	\vspace*{-4pt}
\end{table} 

\section{Group $\mathcal C_i$ ($\bar{1}$)} \label{sec:C_i}

For this point group, we simply adopt the 2nd-order structural tensors given by Zheng \cite{zheng_theory_1994} as follows. 
\begin{align}\label{eq structure_C3i}
	\begin{array}{l}
		\mathbf K_1 =\bm \varepsilon \mathbf i =  \begin{bmatrix}
			0 & 0 & 0 \\
			0 & 0 & 1\\
			0 & -1 & 0
		\end{bmatrix},
		\mathbf K_2 = \bm\varepsilon\mathbf j =  \begin{bmatrix}
			0 & 0 & -1 \\
			0 & 0 & 0\\
			1 & 0 & 0
		\end{bmatrix},
		\mathbf K_3 = \bm\varepsilon\mathbf k  = \begin{bmatrix}
			0 & 1 & 0 \\
			-1 & 0 & 0\\
			0 & 0 & 0
		\end{bmatrix}  
	\end{array}
\end{align}
The Boehler-Liu formulation is used to derive the representations of tensor functions.  

The representation of a 2nd-order symmetric tensor-valued function $\mathbf T (\mathbf C) = \hat{\mathbf T} (\mathbf C, \mathbf K_1, \mathbf K_2, \mathbf K_3 )$ is derived first, where $\mathbf C$ is a 2nd-order symmetric tensor (e.g., Cauchy-Green tensor in continuum mechanics). 
Considering that $\mathbf K_i$ are skew-symmetric, the tensor generators and invariants can be obtained using Table \ref{Table functional basis_3D} and Table \ref{Table generators_3D}. Specifically, the tensor generators are
\begin{align}\label{eq general_T_C_i}
	\begin{array}{l}
		\mathbf I, \mathbf C, \mathbf C^2, \mathbf K_i^2, \mathbf C \mathbf K_i - \mathbf K_i \mathbf C , \ \mathbf C^2 \mathbf K_i - \mathbf K_i \mathbf C^2 , \ \mathbf K_i \mathbf C \mathbf K_i,\ \mathbf K_i \mathbf C \mathbf K_i^2 - \mathbf K_i^2 \mathbf C \mathbf K_i,\\ 
		\mathbf K_1 \mathbf K_2 + \mathbf K_2 \mathbf K_1 ,\ \mathbf K_1 \mathbf K_2^2 - \mathbf K_2^2 \mathbf K_1 ,\ \mathbf K_1^2 \mathbf K_2 - \mathbf K_2 \mathbf K_1^2,\  
		\mathbf K_1 \mathbf K_3 + \mathbf K_3 \mathbf K_1 ,\ \mathbf K_1 \mathbf K_3^2 - \mathbf K_3^2 \mathbf K_1 ,\\ \mathbf K_1^2 \mathbf K_3 - \mathbf K_3 \mathbf K_1^2,\ 
		\mathbf K_2 \mathbf K_3 + \mathbf K_3 \mathbf K_2 ,\ \mathbf K_2 \mathbf K_3^2 - \mathbf K_3^2 \mathbf K_2 ,\ \mathbf K_2^2 \mathbf K_3 - \mathbf K_3 \mathbf K_2^2, \ \text{for}\ i=1,2,3
	\end{array}
\end{align}
The invariants are
\begin{align}\label{eq general_alpha_C_i}
	\begin{array}{l}
		tr \mathbf C, tr \mathbf C^2, tr \mathbf C^3, tr( \mathbf C \mathbf K_i^2), tr( \mathbf C^2 \mathbf K_i^2),tr( \mathbf C^2 \mathbf K_i^2 \mathbf C \mathbf K_i), \\
		tr(\mathbf C \mathbf K_1 \mathbf K_2), tr(\mathbf C \mathbf K_1^2 \mathbf K_2), tr(\mathbf C \mathbf K_1 \mathbf K_2^2), tr(\mathbf C \mathbf K_1 \mathbf K_3), tr(\mathbf C \mathbf K_1^2 \mathbf K_3), \\ tr(\mathbf C \mathbf K_1 \mathbf K_3^2),tr(\mathbf C \mathbf K_2 \mathbf K_3), tr(\mathbf C \mathbf K_2^2 \mathbf K_3), tr(\mathbf C \mathbf K_2 \mathbf K_3^2),
		 \text{\quad for \ } i=1,2,3
	\end{array}
\end{align}
After eliminating the redundant terms in (\ref{eq general_T_C_i}) and (\ref{eq general_alpha_C_i}), the representation of $\hat{\mathbf T}$ is given as 
\begin{align}\label{eq T_Ci}
	\begin{array}{l}
		\hat{\mathbf T} (\mathbf C, \mathbf K_1, \mathbf K_2, \mathbf K_3 ) =  \alpha_1 \mathbf K_1^2 + \alpha_2 \mathbf K_2^2 + \alpha_3 \mathbf K_3^2 + \alpha_{4} (\mathbf K_1 \mathbf K_2 + \mathbf K_2 \mathbf K_1) \\ \quad \quad \quad \quad + \alpha_{5} (\mathbf K_1 \mathbf K_3 + \mathbf K_3 \mathbf K_1) + \alpha_{6} (\mathbf K_2 \mathbf K_3 + \mathbf K_3 \mathbf K_2) 
	\end{array}
\end{align}
and 
\begin{align}\label{eq alpha_Ci}
	\begin{array}{l}
		\alpha_i = \alpha_i (tr( \mathbf C \mathbf K_1^2), tr( \mathbf C \mathbf K_2^2),
		 tr( \mathbf C \mathbf K_3^2), tr( \mathbf C \mathbf K_1 \mathbf K_2), tr( \mathbf C \mathbf K_1 \mathbf K_3), tr(\mathbf C \mathbf K_2 \mathbf K_3) )
	\end{array}
\end{align}

The representation of a scalar-valued function $\hat \psi(\mathbf C, \mathbf K_1, \mathbf K_2, \mathbf K_3)$ follows the same form of (\ref{eq alpha_Ci}).

For this point group, one may use either one structural tensor $\mathbf{P}_2$ or three structural tensors $\mathbf{M}_1,\mathbf{M}_2,\mathbf{M}_3$, as shown in Table \ref{Table structural tensors_3D}. In both cases, the Boehler-Liu formulation is used. We will derive the representations using both approaches below.  

\subsection{Using one structural tensor}
Zheng \cite{zheng_theory_1994} proposed a single structural tensor $\mathbf P_2 = \mathbf i \otimes \mathbf i - \mathbf j \otimes \mathbf j $ for this point group, a 2nd-order symmetric tensor. The representation $\hat{\mathbf T} (\mathbf C, \mathbf P_2 )$ can be obtained using Tables \ref{Table functional basis_3D}-\ref{Table generators_3D}, as
\begin{align}\label{eq general_T_D2h_1}
	\begin{array}{l}
		\hat{\mathbf T} (\mathbf C, \mathbf P_2 ) =\alpha_0 \mathbf I + \alpha_1\mathbf C + \alpha_2 \mathbf C^2 + \alpha_3 \mathbf P_2 + \alpha_4 \mathbf P_2^2 \\
		\quad \quad \ \  +\alpha_5 (\mathbf C \mathbf P_2 +\mathbf P_2 \mathbf C) +\alpha_6 (\mathbf C^2 \mathbf P_2 +\mathbf P_2 \mathbf C^2) +\alpha_7 (\mathbf C \mathbf P_2^2 +\mathbf P_2^2 \mathbf C)
	\end{array}
\end{align}
where $\alpha_i$ is given by
\begin{align}\label{eq alpha_D2h_1}
	\begin{array}{l}
		\alpha_i = \alpha_i (tr \mathbf C, tr \mathbf C^2, tr \mathbf C^3, tr(\mathbf C \mathbf P_2), tr(\mathbf C^2 \mathbf P_2), tr(\mathbf C \mathbf P_2^2), tr(\mathbf C^2 \mathbf P_2^2) )
	\end{array}
\end{align} 

The representation of a scalar-valued function $\hat \psi(\mathbf C, \mathbf P_2)$ follows the same form of (\ref{eq alpha_D2h_1}).

\subsection{Using three structural tensors} \label{Sec 3M_i}
For orthotropic materials, a popular set of structural tensors is $\{ \mathbf{M}_1, \mathbf{M}_2, \mathbf{M}_3 \}$ with $\mathbf M_1 = \mathbf i \otimes \mathbf i$, $\mathbf M_2 = \mathbf j \otimes \mathbf j$ and $\mathbf M_3 = \mathbf k \otimes \mathbf k$. The representation of a tensor-valued function $\hat{\mathbf T}$ was provided by Boehler \cite{boehler_irreducible_1977} as
\begin{align}\label{eq T_D2h_2} 
	\begin{array}{l}
		\hat{\mathbf T}(\mathbf C, \mathbf M_1, \mathbf M_2, \mathbf M_3) = \alpha_1  \mathbf M_1 + \alpha_2  \mathbf M_2 + \alpha_3  \mathbf M_3 + \alpha_4 (\mathbf M_1 \mathbf C + \mathbf C \mathbf M_1) \\ \quad \quad \quad+ \alpha_5 (\mathbf M_2 \mathbf C + \mathbf C \mathbf M_2) + \alpha_6 (\mathbf M_3 \mathbf C + \mathbf C \mathbf M_3) + \alpha_7 \mathbf C^2
	\end{array}
\end{align}
where $\alpha_i $ is
\begin{align} \label{eq alpha_D2h_2}
	\begin{array}{l}
		\alpha_i = \alpha_i(tr(\mathbf C \mathbf M_1), tr(\mathbf C \mathbf M_2), tr(\mathbf C \mathbf M_3), tr(\mathbf C^2 \mathbf M_1), tr(\mathbf C^2 \mathbf M_2), tr(\mathbf C^2 \mathbf M_3), tr \mathbf C^3 )\\ 
		\quad \ =\tilde{\alpha}_i(\mathbf C, \mathbf M_1, \mathbf M_2, \mathbf M_3)
	\end{array}
\end{align}

The representation of a scalar-valued function $\hat \psi$ follows the same form of (\ref{eq alpha_D2h_2}), as
\begin{align} \label{eq psi_D2h_2}
	\begin{array}{l}
		\hat{\psi}(\mathbf C, \mathbf M_1, \mathbf M_2, \mathbf M_3) 
		=  \hat{\psi}(tr(\mathbf C \mathbf M_1), tr(\mathbf C \mathbf M_2), tr(\mathbf C \mathbf M_3), tr(\mathbf C^2 \mathbf M_1), tr(\mathbf C^2 \mathbf M_2), tr(\mathbf C^2 \mathbf M_3), tr \mathbf C^3 )
	\end{array}
\end{align}

These representation formulae for orthotropic materials are very useful for multiple other point groups including $\mathcal{D}_{4h}$, $\mathcal{T}_h$, and $\mathcal{O}_h$ to be introduced below.

\section{Group $\mathcal D_{4h}$ ($4/mmm$)}\label{sec:D_4h}

Zheng \cite{zheng_theory_1994} proposed a 4th-order structural tensor $\mathbb{P}_4$ for this point group. Given the fact that higher-order structural tensors are inconvenient to use, we propose a lower-order structural tensor set $\{\mathbf{M}_1, \mathbf{M}_2, \mathbf{M}_3 \}$ instead, the same as Section \ref{Sec 3M_i}. For this point group, the Man-Goddard reformulation is needed. 

The representations of tensor functions have been provided in Eqs. (\ref{eq T_D2h_2})-(\ref{eq psi_D2h_2}). But additional constraints (\ref{eq_constraint}) must be imposed to the representations. The group generators of this point group are $\mathcal G^* = \{\mathbf C_4,\ \mathbf C_{2x}, \bar{\mathbf I}\}$ in Table A3. The operations $\mathbf C_{2x}$ and $\bar{\mathbf I}$ keep all three structural tensors $\mathbf M_i$ invariant. However, the operation $\mathbf C_4$ permutes $\mathbf M_1$ and $\mathbf M_2$, i.e., $\mathbf C_4 \mathbf M_1 \mathbf C_4^T = \mathbf M_2$ and $\mathbf C_4 \mathbf M_2 \mathbf C_4^T = \mathbf M_1$. Thus, $\mathbf C_4$ would impose additional constraints to the representations based on (\ref{eq_constraint}), as 
\begin{align}
	\begin{array}{l}
		\hat{\mathbf T}(\mathbf C, \mathbf M_1, \mathbf M_2, \mathbf M_3) = \hat{\mathbf T}(\mathbf C, \mathbf M_2, \mathbf M_1, \mathbf M_3), \\
		\hat{\psi}(\mathbf C, \mathbf M_1, \mathbf M_2, \mathbf M_3) = \hat{\psi}(\mathbf C, \mathbf M_2, \mathbf M_1, \mathbf M_3) \\
	\end{array}
\end{align}
Accordingly, the constraints to the coefficients $\tilde{\alpha}_i$ in (\ref{eq alpha_D2h_2}) are 
\begin{align} \label{eq constraints_T_D4h}
	\begin{array}{l}
		\tilde{\alpha}_1(\mathbf C, \mathbf M_1, \mathbf M_2, \mathbf M_3) = \tilde{\alpha}_2(\mathbf C, \mathbf M_2, \mathbf M_1, \mathbf M_3), \\ \tilde{\alpha}_4(\mathbf C, \mathbf M_1, \mathbf M_2, \mathbf M_3) = \tilde{\alpha}_5(\mathbf C, \mathbf M_2, \mathbf M_1, \mathbf M_3), \\
		\tilde{\alpha}_i(\mathbf C, \mathbf M_1, \mathbf M_2, \mathbf M_3) = \tilde{\alpha}_i(\mathbf C, \mathbf M_2, \mathbf M_1, \mathbf M_3) \quad \text{for} \ i=3,6,7 
	\end{array}
\end{align}

\textbf{Remark \ref{sec:D_4h}.1.} Some constraints on the scalar coefficient functions $\tilde{\alpha}_i$ are redundant and should be removed. The reason is that not all constraints are independent. For example, in (\ref{eq constraints_T_D4h}), we have removed a constraint $\tilde{\alpha}_2(\mathbf C, \mathbf M_1, \mathbf M_2, \mathbf M_3) = \tilde{\alpha}_1(\mathbf C, \mathbf M_2, \mathbf M_1, \mathbf M_3) $ because it is equivalent to the first equation in (\ref{eq constraints_T_D4h}). Hence, extra efforts are required to remove redundancy of the representations for each group.

\section{Group $\mathcal T_{h}$ ($m\bar{3}$) } \label{sec:T_h}
This point group has a 4th-order structural tensor $\mathbb T_h$ \cite{zheng_theory_1994}. We adopt the lower-order structural tensor set $\{\mathbf{M}_1, \mathbf{M}_2, \mathbf{M}_3 \}$ in Section \ref{Sec 3M_i}. The representations of tensor functions have been shown in (\ref{eq T_D2h_2})-(\ref{eq psi_D2h_2}). Additional constraints need to be imposed following the Man-Goddard reformulation. 

This point group has four generators $\mathcal G^* = \{\mathbf Q_{p}^{2\pi/3},\ \mathbf C_{2x}, \mathbf C_{2y}, \bar{\mathbf I}\}$ in Table \ref{Table Laue generators}. The operations $\mathbf C_{2x}, \mathbf C_{2y}, \text{ and } \bar{\mathbf I}$ keep all three structural tensors invariant. In contrast, $\mathbf Q_{p}^{2\pi/3}$ permutes them, as
\begin{align}
	\begin{array}{l}
		\mathbf Q^{2\pi/3}_{p} \mathbf M_1 (\mathbf Q^{2\pi/3}_{p})^T = \mathbf M_2, \quad 
		\mathbf Q^{2\pi/3}_{p} \mathbf M_2 (\mathbf Q^{2\pi/3}_{p})^T = \mathbf M_3, \quad 
		\mathbf Q^{2\pi/3}_{p} \mathbf M_3 (\mathbf Q^{2\pi/3}_{p})^T = \mathbf M_1 
	\end{array}
\end{align}
Thus, $\mathbf Q_{p}^{2\pi/3}$ will impose an additional constraint to the representations, as $\hat{\mathbf T}(\mathbf C, \mathbf M_1, \mathbf M_2, \mathbf M_3) = \hat{\mathbf T}(\mathbf C, \mathbf M_2, \mathbf M_3, \mathbf M_1)$. Accordingly, the constraints to the coefficient functions $\tilde{\alpha}_i$ in (\ref{eq T_D2h_2}) are  
\begin{align}\label{eq constraint_T_Th}
	\begin{array}{ll}
		\tilde{\alpha}_1(\mathbf C, \mathbf M_1, \mathbf M_2, \mathbf M_3) = \tilde{\alpha}_3(\mathbf C, \mathbf M_2, \mathbf M_3, \mathbf M_1), & \tilde{\alpha}_2(\mathbf C, \mathbf M_1, \mathbf M_2, \mathbf M_3) = \tilde{\alpha}_1(\mathbf C, \mathbf M_2, \mathbf M_3, \mathbf M_1) \\
		\tilde{\alpha}_4(\mathbf C, \mathbf M_1, \mathbf M_2, \mathbf M_3) = \tilde{\alpha}_6(\mathbf C, \mathbf M_2, \mathbf M_3, \mathbf M_1), &
		\tilde{\alpha}_5(\mathbf C, \mathbf M_1, \mathbf M_2, \mathbf M_3) = \tilde{\alpha}_4(\mathbf C, \mathbf M_2, \mathbf M_3, \mathbf M_1) \\
		\tilde{\alpha}_7(\mathbf C, \mathbf M_1, \mathbf M_2, \mathbf M_3) = \tilde{\alpha}_7(\mathbf C, \mathbf M_2, \mathbf M_3, \mathbf M_1)
	\end{array}
\end{align}

For the scalar-valued function $\hat{\psi}$, the additional constraint on (\ref{eq psi_D2h_2}) requires that $\hat{\psi}(\mathbf C, \mathbf M_1, \mathbf M_2, \mathbf M_3) = \hat{\psi}(\mathbf C, \mathbf M_2, \mathbf M_3, \mathbf M_1)$.  

\section{Group $\mathcal O_h$ ($m\bar{3}m$) } \label{sec:O_h}
Rather than a 4th-order structural tensor $\mathbb O_h$ \cite{zheng_theory_1994}, we adopt the structural tensor set $\{\mathbf M_1, \mathbf M_2, \mathbf M_3 \}$ in Section \ref{Sec 3M_i}. The representations of tensor functions have been shown in (\ref{eq T_D2h_2})-(\ref{eq psi_D2h_2}). Additional constraints need to be imposed following the Man-Goddard reformulation. 

This point group has four group generators $\mathcal G^* = \{\mathbf Q_{p}^{2\pi/3},\ \mathbf C_{4x}, \mathbf C_{2y}, \bar{\mathbf I}\}$ in Table \ref{Table Laue generators}. Among the four group generators, $\mathbf C_{2y} \text{ and } \bar{\mathbf I}$ keep the structural tensors invariant, whereas $\mathbf C_{4x}$ and $\mathbf Q_{p}^{2\pi/3}$ transform them as follows.  
\begin{align}
	\begin{array}{ll}
		\mathbf C_{4x} \mathbf M_1 \mathbf C_{4x}^T = \mathbf M_1, & \mathbf Q^{2\pi/3}_{p} \mathbf M_1 (\mathbf Q^{2\pi/3}_{p})^T = \mathbf M_2\\ 
		\mathbf C_{4x} \mathbf M_2 \mathbf C_{4x}^T = \mathbf M_3, & \mathbf Q^{2\pi/3}_{p} \mathbf M_2 (\mathbf Q^{2\pi/3}_{p})^T = \mathbf M_3\\ 
		\mathbf C_{4x} \mathbf M_3 \mathbf C_{4x}^T = \mathbf M_2, & \mathbf Q^{2\pi/3}_{p} \mathbf M_3 (\mathbf Q^{2\pi/3}_{p})^T = \mathbf M_1\\ 
	\end{array}
\end{align}
Hence, we need to impose additional constraints (\ref{eq_constraint}) for $\mathbf C_{4x}$ and $\mathbf Q_{p}^{2\pi/3}$, respectively. As to the tensor-valued function $\hat{\mathbf T}$, the group generator $\mathbf C_{4x}$ requires that $\hat{\mathbf T}(\mathbf C, \mathbf M_1, \mathbf M_2, \mathbf M_3) = \hat{\mathbf T}(\mathbf C, \mathbf M_1, \mathbf M_3, \mathbf M_2)$, while the group generator $\mathbf Q_{p}^{2\pi/3}$ requires that $\hat{\mathbf T}(\mathbf C, \mathbf M_1, \mathbf M_2, \mathbf M_3) = \hat{\mathbf T}(\mathbf C, \mathbf M_2, \mathbf M_3, \mathbf M_1)$. Accordingly, the additional constraints to the coefficient functions $\tilde{\alpha}_i$ in (\ref{eq T_D2h_2}) are as follows. 
\begin{align}
	\begin{array}{l}
		\text{For the generator } \mathbf{C}_{4x}: \\
		\tilde{\alpha}_3(\mathbf C, \mathbf M_1, \mathbf M_2, \mathbf M_3) = \tilde{\alpha}_2(\mathbf C, \mathbf M_1, \mathbf M_3, \mathbf M_2),\quad
		\tilde{\alpha}_6(\mathbf C, \mathbf M_1, \mathbf M_2, \mathbf M_3) = \tilde{\alpha}_5(\mathbf C, \mathbf M_1, \mathbf M_3, \mathbf M_2), \\
		\tilde{\alpha}_i(\mathbf C, \mathbf M_1, \mathbf M_2, \mathbf M_3) = \tilde{\alpha}_i(\mathbf C, \mathbf M_1, \mathbf M_3, \mathbf M_2) \quad \text{for} \ i=1,4,7 \\
		\text{For the generator } \mathbf{Q}_p^{2\pi/3}: \\
		\tilde{\alpha}_1(\mathbf C, \mathbf M_1, \mathbf M_2, \mathbf M_3) = \tilde{\alpha}_3(\mathbf C, \mathbf M_2, \mathbf M_3, \mathbf M_1), \quad \tilde{\alpha}_2(\mathbf C, \mathbf M_1, \mathbf M_2, \mathbf M_3) = \tilde{\alpha}_1(\mathbf C, \mathbf M_2, \mathbf M_3, \mathbf M_1), \\
		\tilde{\alpha}_4(\mathbf C, \mathbf M_1, \mathbf M_2, \mathbf M_3) = \tilde{\alpha}_6(\mathbf C, \mathbf M_2, \mathbf M_3, \mathbf M_1), \quad
		\tilde{\alpha}_5(\mathbf C, \mathbf M_1, \mathbf M_2, \mathbf M_3) = \tilde{\alpha}_4(\mathbf C, \mathbf M_2, \mathbf M_3, \mathbf M_1),\quad  \\
		\tilde{\alpha}_7(\mathbf C, \mathbf M_1, \mathbf M_2, \mathbf M_3) =
		\tilde{\alpha}_7(\mathbf C, \mathbf M_2, \mathbf M_3, \mathbf M_1)
	\end{array}
\end{align}
Similarly, $\mathbf C_{4x}$ and $\mathbf Q_{p}^{2\pi/3}$ also impose additional constraints to the scalar-valued tensor function $\hat \psi$ in (\ref{eq psi_D2h_2}) as
\begin{align} \label{eq constraint_psi_Oh}
	\hat{\psi}(\mathbf C, \mathbf M_1, \mathbf M_2, \mathbf M_3) = 
	\hat{\psi}(\mathbf C, \mathbf M_1, \mathbf M_3, \mathbf M_2) = \hat{\psi}(\mathbf C, \mathbf M_2, \mathbf M_3, \mathbf M_1)
\end{align} 

\section{Group $\mathcal C_{2h}$ ($2/m$)}\label{sec:C_2h}
For this point group, we simply adopt the two 2nd-order structural tensors $\mathbf P_2 = \mathbf i \otimes \mathbf i - \mathbf j \otimes \mathbf j $ and $\mathbf K_3 = \bm\varepsilon \mathbf k $ proposed by Zheng \cite{zheng_theory_1994}. the Boehler-Liu formulation is used to derive the representations. Considering that $\mathbf P_2$ is symmetric and $\mathbf K_3$ is skew-symmetric, the tensor generators are obtained from Table \ref{Table generators_3D} as
\begin{align}\label{eq general_T_C2h_1}
	\begin{array}{l}
		\mathbf I,\ \mathbf C ,\ \mathbf C^2 ,\ \mathbf P_2 ,\ \mathbf P_2^2 ,\ \mathbf C \mathbf P_2 +\mathbf P_2 \mathbf C,\ 
		\mathbf C^2 \mathbf P_2 +\mathbf P_2 \mathbf C^2,\
		\mathbf C \mathbf P_2^2 +\mathbf P_2^2 \mathbf C ,\\ \mathbf K_3^2 ,\ \mathbf C \mathbf K_3 - \mathbf K_3 \mathbf C,\  \mathbf C^2 \mathbf K_3 - \mathbf K_3 \mathbf C^2
		,\ \mathbf K_3 \mathbf C \mathbf K_3 ,\ \mathbf K_3\mathbf C \mathbf K_3^2 - \mathbf K_3^2\mathbf C \mathbf K_3 ,\\
		\mathbf P_2 \mathbf K_3 - \mathbf K_3 \mathbf P_2,\ \mathbf P_2^2 \mathbf K_3 - \mathbf K_3 \mathbf P_2^2,\  \mathbf K_3 \mathbf P_2 \mathbf K_3,\ \mathbf K_3\mathbf P_2 \mathbf K_3^2 - \mathbf K_3^2\mathbf P_2 \mathbf K_3
	\end{array}
\end{align}
and the invariants are obtained from Table \ref{Table functional basis_3D} as
\begin{align}\label{eq general_alpha_C2h_1}
	\begin{array}{l}
		tr \mathbf C, tr \mathbf C^2, tr \mathbf C^3, tr(\mathbf C \mathbf P_2), tr(\mathbf C^2 \mathbf P_2), tr(\mathbf C \mathbf P_2^2), tr(\mathbf C^2 \mathbf P_2^2),\\
		tr(\mathbf C \mathbf K_3^2), tr(\mathbf C^2 \mathbf K_3^2), tr(\mathbf C^2 \mathbf K_3^2\mathbf C \mathbf K_3), \\
		tr(\mathbf C \mathbf P_2 \mathbf K_3), tr(\mathbf C^2 \mathbf P_2 \mathbf K_3), tr(\mathbf C \mathbf P_2^2 \mathbf K_3), tr(\mathbf C \mathbf K_3^2 \mathbf P_2 \mathbf K_3) 
	\end{array}
\end{align}

There are redundant terms in (\ref{eq general_T_C2h_1}) and (\ref{eq general_alpha_C2h_1}). After eliminating the redundant terms, we obtain 
\begin{align}\label{eq T_C2h_1}
	\begin{array}{l}
		\hat{\mathbf T} (\mathbf C, \mathbf P_2 , \mathbf K_3 ) =\alpha_0 \mathbf I + \alpha_1\mathbf C + \alpha_2 \mathbf C^2 + \alpha_3 \mathbf P_2 + \alpha_4 \mathbf P_2^2 +\alpha_5 (\mathbf C \mathbf P_2 +\mathbf P_2 \mathbf C) \\
		\quad \quad \ \   +\alpha_6 (\mathbf C^2 \mathbf P_2 +\mathbf P_2 \mathbf C^2)
		+\alpha_7 (\mathbf C \mathbf P_2^2 +\mathbf P_2^2 \mathbf C)+ \alpha_{8} (\mathbf C \mathbf K_3 - \mathbf K_3 \mathbf C)  + \alpha_{9} (\mathbf C^2 \mathbf K_3 - \mathbf K_3 \mathbf C^2)\\
		\quad \quad \ \   + \alpha_{10} \mathbf K_3 \mathbf C \mathbf K_3 + \alpha_{11} (\mathbf K_3\mathbf C \mathbf K_3^2 - \mathbf K_3^2\mathbf C \mathbf K_3 ) + \alpha_{12} (\mathbf P_2 \mathbf K_3 - \mathbf K_3 \mathbf P_2)
	\end{array}
\end{align}
and
\begin{align}\label{eq alpha_C2h_1}
	\begin{array}{l}
		\alpha_i = \alpha_i (tr \mathbf C, tr \mathbf C^2, tr(\mathbf C \mathbf P_2), tr(\mathbf C^2 \mathbf P_2), tr(\mathbf C \mathbf P_2^2), tr(\mathbf C \mathbf P_2 \mathbf K_3), tr(\mathbf C^2 \mathbf P_2 \mathbf K_3) )
	\end{array}
\end{align} 

The representation of a scalar-valued function $\hat \psi$ follows the same form of (\ref{eq alpha_C2h_1}).

\section{Group $\mathcal C_{4h}$ ($4/m$)} \label{sec:C_4h}
For this point group, we propose a lower-order structural tensor set $\{\mathbf M_1, \mathbf M_2, \mathbf M_3, \mathbf K_3 \}$, where $\mathbf M_i$ are defined in Section \ref{Sec 3M_i} and $\mathbf K_3 = \bm{\varepsilon} \mathbf k$. The purpose of $\mathbf K_3$ is to break the in-plane reflection symmetry. The Man-Goddard reformulation needs to be used. The group generators are $\mathcal{G}^* = \{\mathbf C_4, \bar{\mathbf I} \}$. The group generator $\bar{\mathbf I}$ keeps all structural tensors invariant. In contrast, the group generator $\mathbf C_4$ keeps $\mathbf M_3$ and $\mathbf K_3$ invariant but permutes $\mathbf M_1 \text{ and } \mathbf M_2$. Hence, $\mathbf C_4$ would impose additional constraints to the representations. 

The tensor generators and invariants can be obtained by adding extra terms related to $\mathbf K_3$ in (\ref{eq T_D2h_2}) and (\ref{eq alpha_D2h_2}). The tensor generators are
\begin{align}\label{eq general_T_C4h}
	\begin{array}{l}
		\mathbf M_i ,\ \mathbf M_i \mathbf C + \mathbf C \mathbf M_i,\ \mathbf C^2 ,\  \mathbf K_3^2 ,\ \mathbf C \mathbf K_3 - \mathbf K_3 \mathbf C,\ \mathbf C^2 \mathbf K_3 - \mathbf K_3 \mathbf C^2,\  \mathbf K_3  \mathbf C \mathbf K_3 ,\ \mathbf K_3 \mathbf C \mathbf K_3^2 -  \mathbf K_3^2 \mathbf C \mathbf K_3, \\ \mathbf M_i\mathbf K_3 - \mathbf K_3 \mathbf M_i ,\ \mathbf M_i^2\mathbf K_3 - \mathbf K_3 \mathbf M_i^2,\ \mathbf K_3 \mathbf M_i \mathbf K_3,\ \mathbf K_3 \mathbf M_i \mathbf K_3^2 - \mathbf K_3^2 \mathbf M_i \mathbf K_3, \text{\ \ for }i=1,2,3 
	\end{array}
\end{align}
and the invariants are
\begin{align} \label{eq general_alpha_C4h}
	\begin{array}{l}
		tr(\mathbf C \mathbf M_i), tr(\mathbf C^2 \mathbf M_i), tr \mathbf C^3, 
		tr(\mathbf C \mathbf K_3^2), tr(\mathbf C^2 \mathbf K_3^2), tr(\mathbf C^2 \mathbf K_3^2 \mathbf C \mathbf K_3),\\ tr(\mathbf C \mathbf M_i \mathbf K_3), tr(\mathbf C^2 \mathbf M_i \mathbf K_3), tr(\mathbf C \mathbf M_i^2 \mathbf K_3), tr(\mathbf C \mathbf K_3^2 \mathbf M_i \mathbf K_3), \text{\ \ for } i=1,2,3
	\end{array}
\end{align}
After eliminating the redundant terms in (\ref{eq general_T_C4h}) and (\ref{eq general_alpha_C4h}), the representation of $\hat{\mathbf T}$ is obtained as
\begin{align}\label{eq T_C4h}
	\begin{array}{l}
		\hat{\mathbf T}(\mathbf C, \mathbf M_1, \mathbf M_2, \mathbf M_3, \mathbf K_3) = \alpha_1  \mathbf M_1 + \alpha_2  \mathbf M_2 + \alpha_3  \mathbf M_3 + \alpha_4 (\mathbf M_1 \mathbf C + \mathbf C \mathbf M_1) \\ \quad \quad \quad+ \alpha_5 (\mathbf M_2 \mathbf C + \mathbf C \mathbf M_2) + \alpha_6 \mathbf C^2 + \alpha_7 ( \mathbf C \mathbf K_3 - \mathbf K_3 \mathbf C)\\  \quad \quad \quad
		+ \alpha_{8} ( \mathbf C^2 \mathbf K_3 - \mathbf K_3 \mathbf C^2) 
		+ \alpha_{9} (\mathbf M_1\mathbf K_3 - \mathbf K_3 \mathbf M_1) + \alpha_{10} (\mathbf M_2\mathbf K_3 - \mathbf K_3 \mathbf M_2)
	\end{array}
\end{align}
where $\alpha_i$ are
\begin{align} \label{eq alpha_C4h}
	\begin{array}{l}
		\alpha_i = \alpha_i(tr(\mathbf C \mathbf M_1), tr(\mathbf C \mathbf M_2), tr(\mathbf C \mathbf M_3), tr(\mathbf C^2 \mathbf M_1), tr(\mathbf C^2 \mathbf M_2),\\ \quad \quad \ 
		tr(\mathbf C \mathbf M_1 \mathbf K_3), tr(\mathbf C^2 \mathbf M_1 \mathbf K_3), tr(\mathbf C \mathbf M_2 \mathbf K_3), tr(\mathbf C^2 \mathbf M_2 \mathbf K_3))\\ 
		\quad \ =\tilde{\alpha}_i(\mathbf C, \mathbf M_1, \mathbf M_2, \mathbf M_3, \mathbf K_3)
	\end{array}
\end{align}

Note that the group generator $\mathbf C_4$ imposes an additional constraint $\hat{\mathbf T}(\mathbf C, \mathbf M_1, \mathbf M_2, \mathbf M_3, \mathbf K_3) = \hat{\mathbf T}(\mathbf C, \mathbf M_2, \mathbf M_1, \mathbf M_3, \mathbf K_3)$ to the representation. Using (\ref{eq T_C4h})-(\ref{eq alpha_C4h}), the additional constraints to coefficient functions are as follows.
\begin{align}
	\begin{array}{l}
		\tilde{\alpha}_1(\mathbf C, \mathbf M_1, \mathbf M_2, \mathbf M_3, \mathbf K_3) = \tilde{\alpha}_2(\mathbf C, \mathbf M_2, \mathbf M_1, \mathbf M_3, \mathbf K_3),\\
		\tilde{\alpha}_4(\mathbf C, \mathbf M_1, \mathbf M_2, \mathbf M_3, \mathbf K_3) = \tilde{\alpha}_5(\mathbf C, \mathbf M_2, \mathbf M_1, \mathbf M_3, \mathbf K_3), \\
		\tilde{\alpha}_{9}(\mathbf C, \mathbf M_1, \mathbf M_2, \mathbf M_3, \mathbf K_3) = \tilde{\alpha}_{10}(\mathbf C, \mathbf M_2, \mathbf M_1, \mathbf M_3, \mathbf K_3),\\
		\tilde{\alpha}_{i}(\mathbf C, \mathbf M_1, \mathbf M_2, \mathbf M_3, \mathbf K_3) = \tilde{\alpha}_{i}(\mathbf C, \mathbf M_2, \mathbf M_1, \mathbf M_3, \mathbf K_3) \quad \text{for} \   i=3,6,7,8
	\end{array}
\end{align}

The representation of a scalar-valued function $\hat \psi$ follows the same form of (\ref{eq alpha_C4h}). 
Moreover, the group generator $\mathbf C_4$ imposes an additional constraint $\hat \psi(\mathbf C, \mathbf M_1, \mathbf M_2, \mathbf M_3, \mathbf K_3) = \hat \psi(\mathbf C, \mathbf M_2, \mathbf M_1, \mathbf M_3, \mathbf K_3)$ to the representation. 

\section{Group $\mathcal C_{3i}$ ($\bar{3}\  \text{or}\  \mathcal S_6$)}\label{sec:C_3i}
Since the structural tensors proposed by Zheng \cite{zheng_theory_1994} include a 4th-order tensor, we need to construct a lower-order structural tensor set for this point group. By defining a vector $\mathbf u = \mathbf i + \mathbf k$ in a high symmetry plane, we propose a structural tensor set $\{\mathbf T_1, \mathbf T_2, \mathbf T_3, \mathbf K_3 \}$ with detailed tensors given as
\begin{align}\label{eq str_c3i}
	\begin{array}{l}
		\mathbf T_1 = \mathbf u \otimes \mathbf u, \quad 
		\mathbf T_2 = \mathbf C_3 \mathbf T_1 \mathbf C_3^T,\quad 
		\mathbf T_3 = \mathbf C_3 \mathbf T_2 \mathbf C_3^T,\quad 
		\mathbf K_3 = \bm\varepsilon \mathbf k
	\end{array}
\end{align}
The purpose of $\mathbf K_3$ is to break the in-plane reflection symmetry. For this point group, the Man-Goddard reformulation is needed to derive the representations. 
This point group has two group generators $\mathcal{G}^* = \{\mathbf C_3, \bar{\mathbf I} \}$. The operation $\bar{\mathbf I}$ keeps all four structural tensors invariant. The operation $\mathbf C_3$ keeps $\mathbf K_3$ invariant but permutes $\mathbf T_1$, $\mathbf T_2$, $\mathbf T_3$ to each other as shown in (\ref{eq str_c3i}). Hence, the operation $\mathbf C_3$ imposes additional constraints to the representations.

Using Tables \ref{Table functional basis_3D}-\ref{Table generators_3D}, the tensor generators are
\begin{align} \label{eq general_T_C3i}
	\begin{array}{l}
		\mathbf I,\ \mathbf{C}, \ \mathbf{C^2},\ \mathbf T_i,\ \mathbf T_i^2,\ \mathbf K_3^2,\ \mathbf C \mathbf T_i + \mathbf T_i \mathbf C,\ \mathbf C^2 \mathbf T_i + \mathbf T_i \mathbf C^2,\ \mathbf C \mathbf T_i^2 + \mathbf T_i^2 \mathbf C,\\  \mathbf T_1 \mathbf T_2 + \mathbf T_2 \mathbf T_1,\ \mathbf T_1^2 \mathbf T_2 + \mathbf T_2 \mathbf T_1^2,\ \mathbf T_1 \mathbf T_2^2 + \mathbf T_2^2 \mathbf T_1,\  \mathbf T_1 \mathbf T_3 + \mathbf T_3 \mathbf T_1,\ \mathbf T_1^2 \mathbf T_3 + \mathbf T_3 \mathbf T_1^2, \\ \mathbf T_1 \mathbf T_3^2 + \mathbf T_3^2 \mathbf T_1,\ \mathbf T_2 \mathbf T_3 + \mathbf T_3 \mathbf T_2,\ \mathbf T_2^2 \mathbf T_3 + \mathbf T_3 \mathbf T_2^2,\ \mathbf T_2 \mathbf T_3^2 + \mathbf T_3^2 \mathbf T_2, \\ \mathbf C \mathbf K_3 -\mathbf K_3 \mathbf C,\ \mathbf C^2 \mathbf K_3 -\mathbf K_3 \mathbf C^2,\ \mathbf K_3  \mathbf C \mathbf K_3,\ \mathbf K_3  \mathbf C \mathbf K_3^2 - \mathbf K_3^2  \mathbf C \mathbf K_3,\\
		\mathbf T_i \mathbf K_3 -\mathbf K_3 \mathbf T_i,\ \mathbf T_i^2 \mathbf K_3 -\mathbf K_3 \mathbf T_i^2,\ \mathbf K_3  \mathbf T_i \mathbf K_3,\ \mathbf K_3  \mathbf T_i \mathbf K_3^2 - \mathbf K_3^2  \mathbf T_i \mathbf K_3, \text{ for  }i=1,2,3 
	\end{array}
\end{align} 
and the invariants are 
\begin{align}\label{eq general_alpha_C3i}
	\begin{array}{l}
		tr\mathbf C, tr \mathbf C^2, tr \mathbf C^3, tr(\mathbf C \mathbf T_i), tr(\mathbf C^2 \mathbf T_i),tr(\mathbf C \mathbf T_i^2), tr(\mathbf C^2 \mathbf T_i^2),
		tr(\mathbf C \mathbf K_3^2),\\ tr(\mathbf C^2 \mathbf K_3^2),tr(\mathbf C^2 \mathbf K_3^2 \mathbf C \mathbf K_3), tr(\mathbf C \mathbf T_1 \mathbf T_2), tr(\mathbf C \mathbf T_1 \mathbf T_3), tr(\mathbf C \mathbf T_2 \mathbf T_3),\\ 
		tr(\mathbf C \mathbf T_i \mathbf K_3), tr(\mathbf C^2 \mathbf T_i \mathbf K_3), tr(\mathbf C \mathbf T_i^2 \mathbf K_3), tr(\mathbf C \mathbf K_3^2 \mathbf T_i \mathbf K_3), \text{ for  }i=1,2,3
	\end{array}
\end{align}
After eliminating the redundant terms in (\ref{eq general_T_C3i}) and (\ref{eq general_alpha_C3i}), the representation of $\hat{\mathbf T}$ is as follows.
\begin{align} \label{eq C3i cons}
	\begin{array}{l}
		\hat{\mathbf T} (\mathbf C, \mathbf T_1, \mathbf T_2, \mathbf T_3, \mathbf K_3) = \alpha_0 \mathbf I + \alpha_1 \mathbf T_1 + \alpha_2 \mathbf T_2 + \alpha_3 \mathbf T_3 + \alpha_4 \mathbf C 
		+ \alpha_5 \mathbf C^2 \\ \quad \quad \quad 
		+ \alpha_{6} (\mathbf C \mathbf T_1 + \mathbf T_1 \mathbf C) + \alpha_{7} (\mathbf C^2 \mathbf T_1 + \mathbf T_1 \mathbf C^2) +  \alpha_{8} (\mathbf C \mathbf T_2 + \mathbf T_2 \mathbf C)\\ 
		\quad \quad \quad 
		+ \alpha_{9} (\mathbf C^2 \mathbf T_2 + \mathbf T_2 \mathbf C^2) +  \alpha_{10} (\mathbf C \mathbf T_3 + \mathbf T_3 \mathbf C) 
		+ \alpha_{11} (\mathbf C^2 \mathbf T_3 + \mathbf T_3 \mathbf C^2)\\ 
		\quad \quad \quad
		+  \alpha_{12} (\mathbf T_1 \mathbf T_2 + \mathbf T_2 \mathbf T_1) 
		+\alpha_{13} (\mathbf T_1 \mathbf T_3 + \mathbf T_3 \mathbf T_1)
		+ \alpha_{14} (\mathbf T_2 \mathbf T_3 + \mathbf T_3 \mathbf T_2) \\ 
		\quad \quad \quad 
		+ \alpha_{15} (\mathbf C \mathbf K_3 -\mathbf K_3 \mathbf C)+ \alpha_{16} (\mathbf C^2 \mathbf K_3 -\mathbf K_3 \mathbf C^2) + \alpha_{17} \mathbf K_3  \mathbf C \mathbf K_3 \\ 
		\quad \quad \quad 
		+ \alpha_{18} ( \mathbf K_3  \mathbf C \mathbf K_3^2 - \mathbf K_3^2  \mathbf C \mathbf K_3)+  \alpha_{19} \mathbf K_3  \mathbf T_1 \mathbf K_3+ \alpha_{20} \mathbf K_3  \mathbf T_2 \mathbf K_3 + \alpha_{21} \mathbf K_3  \mathbf T_3 \mathbf K_3
	\end{array}
\end{align}
where 
\begin{align}\label{eq C3i basis}
	\begin{array}{l}
		\alpha_i = \alpha_i (tr(\mathbf C \mathbf T_1), tr(\mathbf C \mathbf T_2), tr(\mathbf C \mathbf T_3), tr(\mathbf C \mathbf T_1 \mathbf T_2), tr(\mathbf C \mathbf T_1 \mathbf T_3), tr(\mathbf C \mathbf T_2 \mathbf T_3))\\ \ \quad = \tilde{\alpha_i} (\mathbf C, \mathbf T_1, \mathbf T_2, \mathbf T_3, \mathbf K_3)
	\end{array}
\end{align}
Note that the group generator $\mathbf C_3$ imposes an additional constraint $\hat{\mathbf T} (\mathbf C, \mathbf T_1, \mathbf T_2, \mathbf T_3, \mathbf K_3) = \hat{\mathbf T} (\mathbf C, \mathbf T_2, \mathbf T_3, \mathbf T_1, \mathbf K_3)$, which requires the coefficient functions to satisfy
\begin{align}
	\begin{array}{l}
		\tilde \alpha_1(\mathbf C, \mathbf T_1, \mathbf T_2, \mathbf T_3, \mathbf K_3) = \tilde \alpha_3(\mathbf C, \mathbf T_2, \mathbf T_3, \mathbf T_1, \mathbf K_3), \quad
		\tilde \alpha_2(\mathbf C, \mathbf T_1, \mathbf T_2, \mathbf T_3, \mathbf K_3) = \tilde \alpha_1(\mathbf C, \mathbf T_2, \mathbf T_3, \mathbf T_1, \mathbf K_3),\\
		\tilde \alpha_6(\mathbf C, \mathbf T_1, \mathbf T_2, \mathbf T_3, \mathbf K_3) = \tilde \alpha_{10}(\mathbf C, \mathbf T_2, \mathbf T_3, \mathbf T_1, \mathbf K_3), \quad \tilde \alpha_7(\mathbf C, \mathbf T_1, \mathbf T_2, \mathbf T_3, \mathbf K_3) = \tilde \alpha_{11}(\mathbf C, \mathbf T_2, \mathbf T_3, \mathbf T_1, \mathbf K_3)\\
		\tilde \alpha_8(\mathbf C, \mathbf T_1, \mathbf T_2, \mathbf T_3, \mathbf K_3) = \tilde \alpha_{6}(\mathbf C, \mathbf T_2, \mathbf T_3, \mathbf T_1, \mathbf K_3), \quad \tilde \alpha_9(\mathbf C, \mathbf T_1, \mathbf T_2, \mathbf T_3, \mathbf K_3) = \tilde \alpha_{7}(\mathbf C, \mathbf T_2, \mathbf T_3, \mathbf T_1, \mathbf K_3)\\
		\tilde \alpha_{12}(\mathbf C, \mathbf T_1, \mathbf T_2, \mathbf T_3, \mathbf K_3) = \tilde \alpha_{13}(\mathbf C, \mathbf T_2, \mathbf T_3, \mathbf T_1, \mathbf K_3), \quad
		\tilde \alpha_{14}(\mathbf C, \mathbf T_1, \mathbf T_2, \mathbf T_3, \mathbf K_3) = \tilde \alpha_{12}(\mathbf C, \mathbf T_2, \mathbf T_3, \mathbf T_1, \mathbf K_3), \\ \tilde \alpha_{19}(\mathbf C, \mathbf T_1, \mathbf T_2, \mathbf T_3, \mathbf K_3) = \tilde \alpha_{21}(\mathbf C, \mathbf T_2, \mathbf T_3, \mathbf T_1, \mathbf K_3), \quad
		\tilde \alpha_{20}(\mathbf C, \mathbf T_1, \mathbf T_2, \mathbf T_3, \mathbf K_3) = \tilde \alpha_{19}(\mathbf C, \mathbf T_2, \mathbf T_3, \mathbf T_1, \mathbf K_3),\\
		\tilde \alpha_i(\mathbf C, \mathbf T_1, \mathbf T_2, \mathbf T_3, \mathbf K_3) = \tilde \alpha_i(\mathbf C, \mathbf T_2, \mathbf T_3, \mathbf T_1, \mathbf K_3) \quad \text{for} \ i=0,4,5,15,16,17,18
	\end{array}
\end{align} 

The representation of a scalar-valued function $\hat \psi$ follows the same form of (\ref{eq C3i basis}).
In addition, the group generator $\mathbf C_3$ imposes an additional constraint $\hat{\psi} (\mathbf C, \mathbf T_1, \mathbf T_2, \mathbf T_3, \mathbf K_3) = \hat{\psi} (\mathbf C, \mathbf T_2, \mathbf T_3, \mathbf T_1, \mathbf K_3)$ to the representation. 

\section{Group $\mathcal D_{3d}$ ($\bar{3}m$) }\label{sec:D_3d}
Since the structural tensor $\mathbf k \otimes \mathbb{P}_3$ given by Zheng \cite{zheng_theory_1994} is a 4th-order tensor, we propose a lower-order structural tensor set to replace it. Firstly, we define a vector $\mathbf{v} = \mathbf{j+k}$. The structural tensor set $\{\mathbf D_1, \mathbf D_2, \mathbf D_3\}$ is then defined by $\mathbf D_1 = \mathbf v \otimes \mathbf v$, $\mathbf D_2 = \mathbf C_3 \mathbf D_1 \mathbf C_3^T$, and $\mathbf D_3 = \mathbf C_3 \mathbf D_2 \mathbf C_3^T$. In this case, the Man-Goddard reformulation is used.

This point group has three group generators $\mathcal{G}^*=\{\mathbf C_3, \mathbf C_{2x}, \bar{\mathbf I} \}$. The operation $\bar{\mathbf I}$ keeps these structural tensors invariant but $\mathbf C_3$ and $\mathbf C_{2x}$ transform them in the following way. 
\begin{align}\label{eq constraint_D3d}
	\begin{array}{l}
		\mathbf C_3 \mathbf D_1 \mathbf C_3^T = \mathbf D_2,\quad \mathbf C_{2x} \mathbf D_1 \mathbf C_{2x}^T = \mathbf D_1,\\
		\mathbf C_3 \mathbf D_2 \mathbf C_3^T = \mathbf D_3,\quad \mathbf C_{2x} \mathbf D_2 \mathbf C_{2x}^T = \mathbf D_3,\\
		\mathbf C_3 \mathbf D_3 \mathbf C_3^T = \mathbf D_1,\quad \mathbf C_{2x} \mathbf D_3 \mathbf C_{2x}^T = \mathbf D_2
	\end{array}
\end{align}
Therefore, the group generators $\mathbf C_3$ and $\mathbf C_{2x}$ would impose additional constraints to the representations.

All three $\mathbf D_i$ are symmetric. The tensor generators are obtained as
\begin{align} \label{eq general_T_D3d}
	\begin{array}{l}
		\mathbf I,\ \mathbf C,\ \mathbf C^2,\ \mathbf D_i,\ \mathbf D_i^2,\   
		\mathbf C \mathbf D_i + \mathbf D_i \mathbf C, \mathbf C^2 \mathbf D_i + \mathbf D_i \mathbf C^2, \mathbf C \mathbf D_i^2 + \mathbf D_i^2 \mathbf C,\\ 		
		\mathbf D_1 \mathbf D_2 + \mathbf D_2 \mathbf D_1,\ \mathbf D_1^2 \mathbf D_2 + \mathbf D_2 \mathbf D_1^2,\ \mathbf D_1 \mathbf D_2^2 + \mathbf D_2^2 \mathbf D_1,\\ 
		\mathbf D_1 \mathbf D_3 + \mathbf D_3 \mathbf D_1,\ 
		\mathbf D_1^2 \mathbf D_3 + \mathbf D_3 \mathbf D_1^2,\ \mathbf D_1 \mathbf D_3^2 + \mathbf D_3^2 \mathbf D_1,\ \\ 
		\mathbf D_2 \mathbf D_3 + \mathbf D_3 \mathbf D_2,\ \mathbf D_2^2 \mathbf D_3 + \mathbf D_3 \mathbf D_2^2,\ \mathbf D_2 \mathbf D_3^2 + \mathbf D_3^2 \mathbf D_2, \text{\ for \ } i=1,2,3
	\end{array}
\end{align}
and the invariants are
\begin{align}\label{eq general_alpha_D3d}
	\begin{array}{l}
		tr\mathbf C, tr \mathbf C^2, tr \mathbf C^3, tr(\mathbf C \mathbf D_i), tr(\mathbf C^2 \mathbf D_i),tr(\mathbf C \mathbf D_i^2), tr(\mathbf C^2 \mathbf D_i^2), \\
		tr(\mathbf C \mathbf D_1 \mathbf D_2),	tr(\mathbf C \mathbf D_1 \mathbf D_3), tr(\mathbf C \mathbf D_2 \mathbf D_3), \text{\ for \ } i=1,2,3
	\end{array}
\end{align}
After eliminating redundant terms in (\ref{eq general_T_D3d}) and (\ref{eq general_alpha_D3d}), the representation is given as 
\begin{align} \label{eq T_D3d}
	\begin{array}{l}
		\hat{\mathbf T} (\mathbf C, \mathbf D_1, \mathbf D_2, \mathbf D_3) = \alpha_0 \mathbf I + \alpha_1 \mathbf D_1 + \alpha_2 \mathbf D_2 + \alpha_3 \mathbf D_3 + \alpha_4 \mathbf C	 +  \alpha_{5} \mathbf C^2  \\ 
		\quad \quad \quad
		+ \alpha_{6} (\mathbf C \mathbf D_1 + \mathbf D_1 \mathbf C)+ \alpha_{7} (\mathbf C^2 \mathbf D_1 + \mathbf D_1 \mathbf C^2)+ \alpha_{8} (\mathbf C \mathbf D_2 + \mathbf D_2 \mathbf C) \\ 
		\quad \quad \quad
		+ \alpha_{9} (\mathbf C^2 \mathbf D_2 + \mathbf D_2 \mathbf C^2) + \alpha_{10} (\mathbf C \mathbf D_3 + \mathbf D_3 \mathbf C) + \alpha_{11} (\mathbf C^2 \mathbf D_3 + \mathbf D_3 \mathbf C^2)\\ 
		\quad \quad \quad 
		+ \alpha_{12} (\mathbf D_1 \mathbf D_2 + \mathbf D_2 \mathbf D_1)+ \alpha_{13} (\mathbf D_1 \mathbf D_3 + \mathbf D_3 \mathbf D_1) + \alpha_{14} (\mathbf D_2 \mathbf D_3 + \mathbf D_3 \mathbf D_2)
	\end{array}
\end{align}
and 
\begin{align}\label{eq alpha_D3d}
	\begin{array}{l}
		\alpha_i = \alpha_i (tr(\mathbf C \mathbf D_1), tr(\mathbf C \mathbf D_2), tr(\mathbf C \mathbf D_3), tr(\mathbf C \mathbf D_1 \mathbf D_2),	tr(\mathbf C \mathbf D_1 \mathbf D_3), tr(\mathbf C \mathbf D_2 \mathbf D_3)) \\ \quad = \tilde{\alpha}_i (\mathbf C, \mathbf D_1, \mathbf D_2, \mathbf D_3)
	\end{array}
\end{align}

As mentioned earlier, the group generators $\mathbf C_3$ and $\mathbf C_{2x}$ impose additional constraints $\hat{\mathbf T}(\mathbf C, \mathbf D_1, \mathbf D_2, \mathbf D_3) = \hat{\mathbf T}(\mathbf C, \mathbf D_2, \mathbf D_3, \mathbf D_1) $ and $\hat{\mathbf T}(\mathbf C, \mathbf D_1, \mathbf D_2, \mathbf D_3) = \hat{\mathbf T}(\mathbf C, \mathbf D_1, \mathbf D_3, \mathbf D_2) $ to the representation, respectively. The corresponding constraints to the coefficient functions $\tilde{\alpha}_i$ are given as
\begin{align}
	\begin{array}{ll}
		\text{For the generator } \mathbf C_3:\\
		\tilde \alpha_1(\mathbf C, \mathbf D_1, \mathbf D_2, \mathbf D_3) = \tilde \alpha_3(\mathbf C, \mathbf D_2, \mathbf D_3, \mathbf D_1), &
		\tilde \alpha_2(\mathbf C, \mathbf D_1, \mathbf D_2, \mathbf D_3) = \tilde \alpha_1(\mathbf C, \mathbf D_2, \mathbf D_3, \mathbf D_1), \\
		\tilde \alpha_6(\mathbf C, \mathbf D_1, \mathbf D_2, \mathbf D_3) = \tilde \alpha_{10}(\mathbf C, \mathbf D_2, \mathbf D_3, \mathbf D_1), & \tilde \alpha_7(\mathbf C, \mathbf D_1, \mathbf D_2, \mathbf D_3) = \tilde \alpha_{11}(\mathbf C, \mathbf D_2, \mathbf D_3, \mathbf D_1),\\
		\tilde \alpha_8(\mathbf C, \mathbf D_1, \mathbf D_2, \mathbf D_3) = \tilde \alpha_6(\mathbf C, \mathbf D_2, \mathbf D_3, \mathbf D_1), & \tilde \alpha_9(\mathbf C, \mathbf D_1, \mathbf D_2, \mathbf D_3) = \tilde \alpha_7(\mathbf C, \mathbf D_2, \mathbf D_3, \mathbf D_1),\\
		\tilde \alpha_{12}(\mathbf C, \mathbf D_1, \mathbf D_2, \mathbf D_3) = \tilde \alpha_{13}(\mathbf C, \mathbf D_2, \mathbf D_3, \mathbf D_1), & \tilde \alpha_{13}(\mathbf C, \mathbf D_1, \mathbf D_2, \mathbf D_3) = \tilde \alpha_{14}(\mathbf C, \mathbf D_2, \mathbf D_3, \mathbf D_1),\\
		\tilde \alpha_i(\mathbf C, \mathbf D_1, \mathbf D_2, \mathbf D_3) = \tilde \alpha_i(\mathbf C, \mathbf D_2, \mathbf D_3, \mathbf D_1) & \text{for} \ i=0,4,5 \\
		\text{For the generator } \mathbf C_{2x}:\\
		\tilde \alpha_2(\mathbf C, \mathbf D_1, \mathbf D_2, \mathbf D_3) = \tilde \alpha_3(\mathbf C, \mathbf D_1, \mathbf D_3, \mathbf D_2), &
		\tilde \alpha_8(\mathbf C, \mathbf D_1, \mathbf D_2, \mathbf D_3) = \tilde \alpha_{10}(\mathbf C, \mathbf D_1, \mathbf D_3, \mathbf D_2), \\ \tilde \alpha_9(\mathbf C, \mathbf D_1, \mathbf D_2, \mathbf D_3) = \tilde \alpha_{11}(\mathbf C, \mathbf D_1, \mathbf D_3, \mathbf D_2), &
		\tilde \alpha_{12}(\mathbf C, \mathbf D_1, \mathbf D_2, \mathbf D_3) = \tilde \alpha_{13}(\mathbf C, \mathbf D_1, \mathbf D_3, \mathbf D_2), \\
		\tilde \alpha_i(\mathbf C, \mathbf D_1, \mathbf D_2, \mathbf D_3) = \tilde \alpha_i(\mathbf C, \mathbf D_1, \mathbf D_3, \mathbf D_2) & \text{for} \ i=0,1,4,5,6,7,14
	\end{array}
\end{align}

The representation of a scalar-valued function $\hat \psi$ follows the same form of (\ref{eq alpha_D3d}).
Moreover, the additional constraints imposed by the group generators $\mathbf C_3$ and $\mathbf C_{2x}$  are  
\begin{align}\label{eq constraint_psi_D3d}
	\hat \psi(\mathbf C,\mathbf D_1, \mathbf D_2, \mathbf D_3) = \hat \psi(\mathbf C,\mathbf D_2, \mathbf D_3, \mathbf D_1)=\hat \psi(\mathbf C,\mathbf D_1, \mathbf D_3, \mathbf D_2)
\end{align}

\section{Group $\mathcal D_{6h}$ ($6/mmm$) }\label{sec:D_6h}
The structural tensor $\mathbb P_6$ provided by Zheng \cite{zheng_theory_1994} is a 6th-order one. We propose a lower-order structural tensor set $\{\mathbf H_1, \mathbf H_2, \mathbf H_3 \}$ to replace it. The three structural tensors are defined as $\mathbf H_1 = \mathbf i \otimes \mathbf i$, $\mathbf H_2 = \mathbf C_6 \mathbf H_1 \mathbf C_6^T$ and $\mathbf H_3 = \mathbf C_6 \mathbf H_2 \mathbf C_6^T$. In this case, the Man-Goddard reformulation is used. This point group has three group generators $\mathcal G^* = \{\mathbf C_6, \mathbf C_{2x}, \bar{\mathbf I} \}$. The generator $\bar{\mathbf I}$ keeps all three structural tensors invariant. The other two generators $\mathbf C_6$ and $\mathbf C_{2x}$ transform them in the following way.
\begin{align}\label{gen_mult_str_D6h}
	\begin{array}{ll}
		\mathbf C_6 \mathbf H_1 \mathbf C_6^T = \mathbf H_2,& \mathbf C_{2x} \mathbf H_1 \mathbf C_{2x}^T = \mathbf H_1,\\
		\mathbf C_6 \mathbf H_2 \mathbf C_6^T = \mathbf H_3,& \mathbf C_{2x} \mathbf H_2 \mathbf C_{2x}^T = \mathbf H_3,\\
		\mathbf C_6 \mathbf H_3 \mathbf C_6^T = \mathbf H_1,& \mathbf C_{2x} \mathbf H_3 \mathbf C_{2x}^T = \mathbf H_2
	\end{array}
\end{align} 
As it is obvious from (\ref{gen_mult_str_D6h}), these two group generators would impose additional constraints to the representations. 

All three structural tensors $\mathbf H_i$ where $i= 1,2,3$ are symmetric. The tensor generators are given as
\begin{align} \label{eq general_T_D6h}
	\begin{array}{l}
		\mathbf I,\ \mathbf C,\ \mathbf C^2,\ \mathbf H_i,\ \mathbf H_i^2,\   
		\mathbf C \mathbf H_i + \mathbf H_i \mathbf C, \mathbf C^2 \mathbf H_i + \mathbf H_i \mathbf C^2, \mathbf C \mathbf H_i^2 + \mathbf H_i^2 \mathbf C,\\ 		
		\mathbf H_1 \mathbf H_2 + \mathbf H_2 \mathbf H_1,\ \mathbf H_1^2 \mathbf H_2 + \mathbf H_2 \mathbf H_1^2,\ \mathbf H_1 \mathbf H_2^2 + \mathbf H_2^2 \mathbf H_1,\\ 
		\mathbf H_1 \mathbf H_3 + \mathbf H_3 \mathbf H_1,\ 
		\mathbf H_1^2 \mathbf H_3 + \mathbf H_3 \mathbf H_1^2,\ \mathbf H_1 \mathbf H_3^2 + \mathbf H_3^2 \mathbf H_1,\ \\ 
		\mathbf H_2 \mathbf H_3 + \mathbf H_3 \mathbf H_2,\ \mathbf H_2^2 \mathbf H_3 + \mathbf H_3 \mathbf H_2^2,\ \mathbf H_2 \mathbf H_3^2 + \mathbf H_3^2 \mathbf H_2, \text{\ for \ } i=1,2,3
	\end{array}
\end{align} 
and the invariants are
\begin{align}\label{eq general_alpha_D6h}
	\begin{array}{l}
		tr\mathbf C, tr \mathbf C^2, tr \mathbf C^3, tr(\mathbf C \mathbf H_i), tr(\mathbf C^2 \mathbf H_i),tr(\mathbf C \mathbf H_i^2), tr(\mathbf C^2 \mathbf H_i^2), \\
		tr(\mathbf C \mathbf H_1 \mathbf H_2),	tr(\mathbf C \mathbf H_1 \mathbf H_3), tr(\mathbf C \mathbf H_2 \mathbf H_3), \text{\ for \ } i=1,2,3
	\end{array}
\end{align}

After eliminating the redundant terms in (\ref{eq general_T_D6h}) and (\ref{eq general_alpha_D6h}), the representation of $\hat{\mathbf T}$ is given as
\begin{align} \label{eq D6h cons}
	\begin{array}{l}
		\hat{\mathbf T} (\mathbf C, \mathbf H_1, \mathbf H_2, \mathbf H_3) =  \alpha_0 \mathbf I + \alpha_1 \mathbf H_1 + \alpha_2 \mathbf H_2 + \alpha_3 \mathbf H_3 +  \alpha_4 \mathbf C +  \alpha_5 \mathbf C^2\\ \quad \quad \quad 
		+ \alpha_{6} (\mathbf C \mathbf H_1 + \mathbf H_1 \mathbf C) + \alpha_{7} (\mathbf C^2 \mathbf H_1 + \mathbf H_1 \mathbf C^2) +  \alpha_{8} (\mathbf C \mathbf H_2 + \mathbf H_2 \mathbf C)\\ 
		\quad \quad \quad 
		+ \alpha_{9} (\mathbf C^2 \mathbf H_2 + \mathbf H_2 \mathbf C^2) +  \alpha_{10} (\mathbf C \mathbf H_3 + \mathbf H_3 \mathbf C)
		+ \alpha_{11} (\mathbf C^2 \mathbf H_3 + \mathbf H_3 \mathbf C^2)
	\end{array}
\end{align} 
and 
\begin{align}\label{eq D6h basis}
	\begin{array}{l}
		\alpha_i = \alpha_i (tr \mathbf C, tr \mathbf C^2, tr \mathbf C^3, tr(\mathbf C \mathbf H_1), tr(\mathbf C^2 \mathbf H_1), tr(\mathbf C \mathbf H_2), tr(\mathbf C^2 \mathbf H_2), tr(\mathbf C \mathbf H_3), tr(\mathbf C^2 \mathbf H_3))\\
		\quad \quad 
		=\tilde \alpha_i(\mathbf C, \mathbf H_1, \mathbf H_2, \mathbf H_3)
	\end{array}
\end{align}
As mentioned earlier, additional constraints are required to be imposed. The group generators $\mathbf C_6$ and $\mathbf C_{2x}$ require that $\hat{\mathbf T}(\mathbf C, \mathbf H_1, \mathbf H_2, \mathbf H_3) = \hat{\mathbf T}(\mathbf C, \mathbf H_2, \mathbf H_3, \mathbf H_1) $ and $\hat{\mathbf T}(\mathbf C, \mathbf H_1, \mathbf H_2, \mathbf H_3) = \hat{\mathbf T}(\mathbf C, \mathbf H_1, \mathbf H_3, \mathbf H_2) $, respectively. Consequently, we can find the constraints to the coefficient functions $\tilde{\alpha}_i$ as 
\begin{align}
	\begin{array}{ll}
		\text{For the generator } \mathbf C_{6}:\\
		\tilde \alpha_1(\mathbf C, \mathbf H_1, \mathbf H_2, \mathbf H_3) = \tilde \alpha_3(\mathbf C, \mathbf H_2, \mathbf H_3, \mathbf H_1), & \tilde \alpha_2(\mathbf C, \mathbf H_1, \mathbf H_2, \mathbf H_3) = \tilde \alpha_1(\mathbf C, \mathbf H_2, \mathbf H_3, \mathbf H_1),\\
		\tilde \alpha_6(\mathbf C, \mathbf H_1, \mathbf H_2, \mathbf H_3) = \tilde \alpha_{10}(\mathbf C, \mathbf H_2, \mathbf H_3, \mathbf H_1), & \tilde \alpha_8(\mathbf C, \mathbf H_1, \mathbf H_2, \mathbf H_3) = \tilde \alpha_{6}(\mathbf C, \mathbf H_2, \mathbf H_3, \mathbf H_1),\\
		\tilde \alpha_7(\mathbf C, \mathbf H_1, \mathbf H_2, \mathbf H_3) = \tilde \alpha_{11}(\mathbf C, \mathbf H_2, \mathbf H_3, \mathbf H_1), & \tilde \alpha_9(\mathbf C, \mathbf H_1, \mathbf H_2, \mathbf H_3) = \tilde \alpha_{7}(\mathbf C, \mathbf H_2, \mathbf H_3, \mathbf H_1),\\
		\tilde \alpha_i(\mathbf C, \mathbf H_1, \mathbf H_2, \mathbf H_3) = \tilde \alpha_i(\mathbf C, \mathbf H_2, \mathbf H_3, \mathbf H_1) & \text{for} \ i=0,4,5 \\
		\text{For the generator } \mathbf C_{2x}:\\
		\tilde \alpha_2(\mathbf C, \mathbf H_1, \mathbf H_2, \mathbf H_3) = \tilde \alpha_3(\mathbf C, \mathbf H_1, \mathbf H_3, \mathbf H_2), & \tilde \alpha_8(\mathbf C, \mathbf H_1, \mathbf H_2, \mathbf H_3) = \tilde \alpha_{10}(\mathbf C, \mathbf H_1, \mathbf H_3, \mathbf H_2),\\
		\tilde \alpha_9(\mathbf C, \mathbf H_1, \mathbf H_2, \mathbf H_3) = \tilde \alpha_{11}(\mathbf C, \mathbf H_1, \mathbf H_3, \mathbf H_2), & \\
		\tilde \alpha_i(\mathbf C, \mathbf H_1, \mathbf H_2, \mathbf H_3) = \tilde \alpha_i(\mathbf C, \mathbf H_1, \mathbf H_3, \mathbf H_2) & \text{for} \ i=0,1,4,5,6,7 \\
	\end{array}
\end{align}

The representation of a scalar-valued function $\hat \psi$ follows the same form of (\ref{eq D6h basis}).
Moreover, the additional constraints imposed by the group generators $\mathbf C_6$ and $\mathbf C_{2x}$ are
\begin{align}\label{eq psi_D6h_constraint}
	\hat \psi(\mathbf C,\mathbf H_1, \mathbf H_2, \mathbf H_3) = \hat \psi(\mathbf C,\mathbf H_2, \mathbf H_3, \mathbf H_1) = \hat \psi(\mathbf C,\mathbf H_1, \mathbf H_3, \mathbf H_2)
\end{align}

\section{Group $\mathcal C_{6h}$ ($6/m$) }\label{sec:C_6h}
The structural tensors provided by Zheng \cite{zheng_theory_1994} involve a 6th-order tensor. We propose a lower-order structural tensor set $\{\mathbf H_1, \mathbf H_2, \mathbf H_3, \mathbf K_3 \}$ for this point group. Herein, $\mathbf H_i$ are identical to that in Section \ref{sec:D_6h} and $\mathbf K_3 = \bm\varepsilon \mathbf k$ is introduced to break the in-plane reflection symmetry. In this case, the Man-Goddard reformulation is used. There are only two group generators $\mathcal G^* = \{\mathbf C_6, \bar{\mathbf I} \}$. The generator $\bar{\mathbf I}$ keeps all structural tensors invariant; whereas $\mathbf C_6$ transform them in the following way.
\begin{align}\label{gen_mult_str_C6h}
	\begin{array}{l} 
		\mathbf C_6 \mathbf H_1 \mathbf C_6^T = \mathbf H_2, \quad
		\mathbf C_6 \mathbf H_2 \mathbf C_6^T = \mathbf H_3, \quad
		\mathbf C_6 \mathbf H_3 \mathbf C_6^T = \mathbf H_1, \quad
		\mathbf C_6 \mathbf K_3 \mathbf C_6^T = \mathbf K_3
	\end{array}
\end{align} 
As it is obvious from (\ref{gen_mult_str_C6h}), $\mathbf C_6$ keeps $\mathbf K_3 $ invariant but permutes $\mathbf H_1, \mathbf H_2 \text{ and } \mathbf H_3$. Hence, the generator $\mathbf C_6$ would impose additional constraints to the representations. 

The representation of a tensor-valued function $\hat{\mathbf T}$ is considered first. We can start with the representations (\ref{eq D6h cons}) and (\ref{eq D6h basis}) for the point group $\mathcal D_{6h}$ and add additional terms related to $\mathbf K_3$. The tensor generators are given as 
\begin{align} \label{eq general_T_C6h}
	\begin{array}{l}
		\mathbf I,\ \mathbf C,\ \mathbf C^2,\ \mathbf H_i,\ 
		\mathbf C \mathbf H_i + \mathbf H_i \mathbf C, \ \mathbf C^2 \mathbf H_i + \mathbf H_i \mathbf C^2,\ \mathbf K_3^2 , \\ \mathbf C \mathbf K_3 - \mathbf K_3 \mathbf C,\
		\mathbf C^2 \mathbf K_3 - \mathbf K_3 \mathbf C^2,\ \mathbf K_3 \mathbf C  \mathbf K_3 ,\ \mathbf K_3 \mathbf C  \mathbf K_3^2 -\mathbf K_3^2 \mathbf C \mathbf K_3, \\ 
		\mathbf H_i \mathbf K_3 - \mathbf K_3 \mathbf H_i,\ \mathbf H_i^2 \mathbf K_3 - \mathbf K_3 \mathbf H_i^2,\ \mathbf K_3 \mathbf H_i  \mathbf K_3,\ \mathbf K_3 \mathbf H_i  \mathbf K_3^2 -\mathbf K_3^2 \mathbf H_i \mathbf K_3, \text{\ for \ } i=1,2,3
	\end{array}
\end{align} 
and the invariants are
\begin{align}\label{eq general_alpha_C6h}
	\begin{array}{l}
		tr \mathbf C, tr \mathbf C^2, tr \mathbf C^3, tr(\mathbf C \mathbf H_i), tr(\mathbf C^2 \mathbf H_i), tr(\mathbf C \mathbf K_3^2),tr(\mathbf C^2 \mathbf K_3^2), tr(\mathbf C^2 \mathbf K_3^2 \mathbf C \mathbf K_3),\\ 
		tr(\mathbf C \mathbf H_i \mathbf K_3),tr(\mathbf C^2 \mathbf H_i \mathbf K_3),tr(\mathbf C \mathbf H_i^2 \mathbf K_3), tr(\mathbf C \mathbf K_3^2 \mathbf H_i \mathbf K_3), \text{\ for \ } i=1,2,3
	\end{array}
\end{align}
After eliminating redundant terms in (\ref{eq general_T_C6h}) and (\ref{eq general_alpha_C6h}), the representation of $\hat{\mathbf T}$ is expressed as 
\begin{align} \label{eq T_C6h}
	\begin{array}{l}
		\hat{\mathbf T} (\mathbf C, \mathbf H_1, \mathbf H_2, \mathbf H_3, \mathbf K_3) =  \alpha_0 \mathbf I + \alpha_1 \mathbf H_1 + \alpha_2 \mathbf H_2 + \alpha_3 \mathbf H_3 +  \alpha_4 \mathbf C +  \alpha_5 \mathbf C^2\\ \quad \quad \quad 
		+ \alpha_{6} (\mathbf C \mathbf H_1 + \mathbf H_1 \mathbf C) + \alpha_{7} (\mathbf C^2 \mathbf H_1 + \mathbf H_1 \mathbf C^2) +  \alpha_{8} (\mathbf C \mathbf H_2 + \mathbf H_2 \mathbf C)
		+ \alpha_{9} (\mathbf C^2 \mathbf H_2 + \mathbf H_2 \mathbf C^2)\\ 
		\quad \quad \quad 
		+  \alpha_{10} (\mathbf C \mathbf H_3 + \mathbf H_3 \mathbf C)
		+ \alpha_{11} (\mathbf C^2 \mathbf H_3 + \mathbf H_3 \mathbf C^2)+ \alpha_{12} (\mathbf C \mathbf K_3 - \mathbf K_3 \mathbf C) \\ \quad \quad \quad
		+ \alpha_{13} (\mathbf C^2 \mathbf K_3 - \mathbf K_3 \mathbf C^2) + \alpha_{14} (\mathbf K_3 \mathbf C  \mathbf K_3 ) + \alpha_{15} (\mathbf K_3 \mathbf C  \mathbf K_3^2 -\mathbf K_3^2 \mathbf C \mathbf K_3) 
	\end{array}
\end{align} 
and 
\begin{align}\label{eq alpha_C6h}
	\begin{array}{l}
		\alpha_i = \alpha_i (tr \mathbf C, tr \mathbf C^2, tr \mathbf C^3, tr(\mathbf C \mathbf H_1), tr(\mathbf C^2 \mathbf H_1), tr(\mathbf C \mathbf H_2), tr(\mathbf C^2 \mathbf H_2), tr(\mathbf C \mathbf H_3), tr(\mathbf C^2 \mathbf H_3)) \\ \quad
		=\tilde \alpha_i(\mathbf C, \mathbf H_1, \mathbf H_2, \mathbf H_3, \mathbf K_3)
	\end{array}
\end{align}

As mentioned earlier, the group generator $\mathbf C_6$ imposes an additional constraint $\hat{\mathbf T}(\mathbf C, \mathbf H_1, \mathbf H_2, \mathbf H_3, \mathbf K_3) = \hat{\mathbf T}(\mathbf C, \mathbf H_2, \mathbf H_3, \mathbf H_1, \mathbf K_3)$ to the representation, which requires the coefficient functions to satisfy the following constraints.
\begin{align}\label{eq T_constraint_C6h}
	\begin{array}{l}
		\tilde \alpha_1(\mathbf C, \mathbf H_1, \mathbf H_2, \mathbf H_3, \mathbf K_3) = \tilde \alpha_3(\mathbf C, \mathbf H_2, \mathbf H_3, \mathbf H_1, \mathbf K_3),\\
		\tilde \alpha_2(\mathbf C, \mathbf H_1, \mathbf H_2, \mathbf H_3, \mathbf K_3) = \tilde \alpha_1(\mathbf C, \mathbf H_2, \mathbf H_3, \mathbf H_1, \mathbf K_3),\\
		\tilde \alpha_6(\mathbf C, \mathbf H_1, \mathbf H_2, \mathbf H_3, \mathbf K_3) = \tilde \alpha_{10}(\mathbf C, \mathbf H_2, \mathbf H_3, \mathbf H_1, \mathbf K_3),\\
		\tilde \alpha_8(\mathbf C, \mathbf H_1, \mathbf H_2, \mathbf H_3, \mathbf K_3) = \tilde \alpha_{6}(\mathbf C, \mathbf H_2, \mathbf H_3, \mathbf H_1, \mathbf K_3),\\
		\tilde \alpha_7(\mathbf C, \mathbf H_1, \mathbf H_2, \mathbf H_3, \mathbf K_3) = \tilde \alpha_{11}(\mathbf C, \mathbf H_2, \mathbf H_3, \mathbf H_1, \mathbf K_3),\\
		\tilde \alpha_9(\mathbf C, \mathbf H_1, \mathbf H_2, \mathbf H_3, \mathbf K_3) = \tilde \alpha_{7}(\mathbf C, \mathbf H_2, \mathbf H_3, \mathbf H_1, \mathbf K_3),\\
		\tilde \alpha_i(\mathbf C, \mathbf H_1, \mathbf H_2, \mathbf H_3, \mathbf K_3) = \tilde \alpha_i(\mathbf C, \mathbf H_2, \mathbf H_3, \mathbf H_1, \mathbf K_3) \quad \text{for} \ i=0,4,5,12,13,14,15\\
	\end{array}
\end{align}

The representation of a scalar-valued function $\hat \psi$ follows the same form of (\ref{eq alpha_C6h}). Moreover, the additional constraint imposed by $\mathbf C_6$ is
\begin{align}\label{eq psi_C6h_constraint}
	\hat \psi(\mathbf C,\mathbf H_1, \mathbf H_2, \mathbf H_3, \mathbf K_3) = \hat \psi(\mathbf C,\mathbf H_2, \mathbf H_3, \mathbf H_1, \mathbf K_3)
\end{align}

\section{Continuous groups}\label{sec:Continuous}
In this section, we provide the representations of tensor functions for three centrosymmetric continuous groups $\mathcal C_{\infty h}, \mathcal D_{\infty h} $ (transversely isotropic) and $\mathcal K_{\infty}$ (isotropic). For simplicity purposes, we only present the final results. Most of the results can be found in the literature. 

For the transversely isotropic group $\mathcal C_{\infty h}$, Zheng proposed $\mathbf K_3 = \bm\varepsilon \mathbf k$ as the structural tensor. The representation of $\hat{\mathbf T}$ is 
\begin{align} \label{eq T_Cinfh}
	\begin{array}{l}
		\hat{\mathbf T} (\mathbf C,\mathbf K_3) =  \alpha_0 \mathbf I + \alpha_1 \mathbf C + \alpha_2 \mathbf C^2 + \alpha_3 \mathbf K_3^2 
		+ \alpha_{4} (\mathbf C \mathbf K_3 - \mathbf K_3 \mathbf C) \\ \quad \quad \quad
		+ \alpha_{5} (\mathbf C^2 \mathbf K_3 - \mathbf K_3 \mathbf C^2) + \alpha_{6} (\mathbf K_3 \mathbf C  \mathbf K_3 ) + \alpha_{7} (\mathbf K_3 \mathbf C  \mathbf K_3^2 -\mathbf K_3^2 \mathbf C \mathbf K_3) 
	\end{array}
\end{align} 
and 
\begin{align}\label{eq alpha_Cinfh}
	\begin{array}{l}
		\alpha_i = \alpha_i (tr \mathbf C, tr \mathbf C^2, tr \mathbf C^3, tr(\mathbf C \mathbf K_3^2), tr(\mathbf C^2 \mathbf K_3^2), tr(\mathbf C^2 \mathbf K_3^2 \mathbf C \mathbf K_3))
	\end{array}
\end{align}
The representation of a scalar-valued function $\hat \psi$ is the same as (\ref{eq alpha_Cinfh}).

For the transversely isotropic group $\mathcal D_{\infty h}$, Boehler \cite{Boehler_applications_1987} proposed $\mathbf M_3 = \mathbf k \otimes \mathbf k$ as the structural tensor. The representation of a tensor-valued function $\hat{\mathbf T}$ is 
\begin{align} \label{eq T_Dinfh}
	\begin{array}{l}
		\hat{\mathbf T} (\mathbf C,\mathbf M_3) =  \alpha_0 \mathbf I + \alpha_1 \mathbf C + \alpha_2 \mathbf C^2 + \alpha_3 \mathbf M_3
		+ \alpha_{4} (\mathbf C \mathbf M_3 + \mathbf M_3 \mathbf C) 
		+ \alpha_{5} (\mathbf C^2 \mathbf M_3 + \mathbf M_3 \mathbf C^2)
	\end{array}
\end{align} 
and 
\begin{align}\label{eq alpha_Dinfh}
	\begin{array}{l}
		\alpha_i = \alpha_i (tr \mathbf C, tr \mathbf C^2, tr \mathbf C^3, tr(\mathbf C \mathbf M_3), tr(\mathbf C^2 \mathbf M_3))
	\end{array}
\end{align}
The representation of a scalar-valued function $\hat{\psi}$ is the same as (\ref{eq alpha_Dinfh}).

Finally, for the isotropic group $\mathcal K_{\infty}$, the representations are well known as 
\begin{align} \label{eq T_Kinf}
	\begin{array}{l}
		\mathbf T (\mathbf C ) =  \alpha_0 \mathbf I + \alpha_1 \mathbf C + \alpha_2 \mathbf C^2 
	\end{array}
\end{align} 
and 
\begin{align}\label{eq alpha_Kinf}
	\begin{array}{l}
		\alpha_i = \alpha_i (tr \mathbf C, tr \mathbf C^2, tr \mathbf C^3)
	\end{array}
\end{align}
The representation of a scalar-valued function $\psi$ is the same as (\ref{eq alpha_Kinf}).

\section{Conclusion}
\sloppy In this work, we present a systematic study on the representation of tensor functions using lower-order structural tensor sets for 3D centrosymmetric point groups. 
The traditional representation theory by Boehler and Liu involves higher-order structural tensors that are inconvenient to use for constitutive modeling of anisotropic materials. Based on a reformulated representation theory by Man and Goddard, we propose lower-order structural tensor sets for 3D centrosymmetric point groups and derive the representations of scalar- and 2nd-order symmetric tensor-valued functions for each group. Among the 14 centrosymmetric groups in 3D space, six groups ($\mathcal{C}_i, \mathcal{C}_{2h}, \mathcal{D}_{2h}, \mathcal{C}_{\infty h}, \mathcal{D}_{\infty h}, \text{ and } \mathcal{K}_{h}$) have lower-order structural tensors so the original Boehler-Liu formulation is used. In contrast, for the eight groups ($\mathcal{C}_{4h}, \mathcal{D}_{4h}, \mathcal{C}_{3i}, \mathcal{D}_{3d}, \mathcal{C}_{6h}, \mathcal{D}_{6h}, \mathcal{T}_{h}, \text{ and } \mathcal{O}_{h}$) involving higher-order structural tensors, the Man-Goddard reformulation and our proposed lower-order structural tensor sets should be used. The key difference between the Boehler-Liu formulation and Man-Goddard reformulation is that the latter relaxes symmetry constraints to structural tensors but requires additional constraints to the representations afterwards. The representation theory developed in this work provides explicit expressions of tensor functions for constitutive modeling of anisotropic materials. For scalar-valued and 2nd-order symmetric tensor-valued functions, the presented theory is applicable to all 3D point groups because their representations are determined by the corresponding centrosymmetric groups. Certainly, the structural tensor sets are non-unique. Researchers can devise new structural tensor sets and derive the representations following a similar procedure. Future research can be towards developing specific constitutive laws for anisotropic materials and integrating the presented theory with artificial intelligence to enable data-driven constitutive modeling. 

\section*{Acknowledgments}
This work was supported by the National Science Foundation (NSF) under Grant No. CMMI-2244952. The authors gratefully acknowledge Professor Chi-Sing Man from University of Kentucky for valuable input to this work.

\section*{Declaration of Interest Statement}
The authors declare that they have no known competing financial interests or personal relationships that could have appeared to influence the work reported in this paper.

\bibliographystyle{unsrt}
\bibliography{Anisotropic_3D.bib}
\newpage{} 
\appendix
\section*{Appendix} \label{appendix:A}
\renewcommand{\thetable}{A\arabic{table}}
\setcounter{table}{0} % Reset table counter

\begin{table}[h]
	\caption{Invariants in the 3D isotropic irreducible functional bases of $\mathbf A_i$ and  $\mathbf W_i$ \cite{zheng_theory_1994, smith_isotropic_1971}.}
	\label{Table functional basis_3D}
	\centering
	\begin{tabular}{|c|c|}%%%The number of columns has to be defined here
		\hline
		Variables & Invariants\\
		\hline
		$\mathbf A$ & $tr\mathbf A$,$tr\mathbf A^2$,$tr\mathbf A^3$\\
		$\mathbf A_1$,$\mathbf A_2$ & $tr(\mathbf A_1 \mathbf A_2)$, $tr(\mathbf A_1^2 \mathbf A_2)$,  $tr(\mathbf A_1 \mathbf A_2^2)$,  $tr(\mathbf A_1^2 \mathbf A_2^2)$\\
		$\mathbf A_1$,$\mathbf A_2$,$\mathbf A_3$ & $tr(\mathbf A_1 \mathbf A_2 \mathbf A_3)$\\
		$\mathbf W$ & $tr\mathbf W^2$\\
		$\mathbf A$,$\mathbf W$ & $tr(\mathbf A \mathbf W^2)$, $tr(\mathbf A^2 \mathbf W^2)$,  $tr(\mathbf A^2 \mathbf W^2 \mathbf A \mathbf W)$\\
		$\mathbf A_1$,$\mathbf A_2$,$\mathbf W$ & $tr(\mathbf A_1 \mathbf A_2 \mathbf W)$ , $tr(\mathbf A_1^2 \mathbf A_2 \mathbf W)$, $tr(\mathbf A_1 \mathbf A_2^2 \mathbf W)$, $tr(\mathbf A_1 \mathbf W^2 \mathbf A_2 \mathbf W)$\\
		$\mathbf W_1$,$\mathbf W_2$ & $tr(\mathbf W_1 \mathbf W_2)$\\
		$\mathbf A$,$\mathbf W_1$,$\mathbf W_2$ & $tr(\mathbf{A} \mathbf{W_1} \mathbf{W_2})$, $tr(\mathbf{A} \mathbf{W_1^2} \mathbf{W_2})$, $tr(\mathbf{A} \mathbf{W_1} \mathbf{W_2^2})$ \\
		$\mathbf W_1$,$\mathbf W_2$,$\mathbf W_3$ & $tr(\mathbf W_1 \mathbf W_2 \mathbf W_3)$ \\
		\hline
	\end{tabular}
	\vspace*{-4pt}
\end{table}  \FloatBarrier

\begin{table}[h]
	\caption{Tensor generators in the 3D irreducible representations for isotropic 2nd-order symmetric tensor-valued functions of $\mathbf A_i$ and $\mathbf W_i$ \cite{zheng_theory_1994, smith_isotropic_1971}.}
	\label{Table generators_3D}
	\centering
	\begin{tabular}{|c|c|}%%%The number of columns has to be defined here
		\hline
		Variables & Generators\\
		\hline
		& $\mathbf I$\\
		$\mathbf A$ & $\mathbf A$, $\mathbf A^2$\\
		$\mathbf W$ & $\mathbf W^2$\\
		$\mathbf A_1$,$\mathbf A_2$ & $\mathbf A_1 \mathbf A_2 + \mathbf A_2 \mathbf A_1 $, $\mathbf A_1^2 \mathbf A_2 + \mathbf A_2 \mathbf A_1^2 $,  $\mathbf A_1 \mathbf A_2^2 + \mathbf A_2^2 \mathbf A_1 $\\
		$\mathbf A$,$\mathbf W$ & $\mathbf A \mathbf W - \mathbf W \mathbf A$, $\mathbf A^2 \mathbf W - \mathbf W \mathbf A^2$, $\mathbf W \mathbf A \mathbf W$, $\mathbf W \mathbf A \mathbf W^2 - \mathbf W^2 \mathbf A \mathbf W$ \\
		$\mathbf W_1$,$\mathbf W_2$ & $\mathbf W_1 \mathbf W_2 + \mathbf W_2 \mathbf W_1 $, $\mathbf W_1 \mathbf W_2^2 - \mathbf W_2^2 \mathbf W_1$, $\mathbf W_1^2 \mathbf W_2 - \mathbf W_2 \mathbf W_1^2$\\
		\hline
	\end{tabular}
	\vspace*{-4pt}
\end{table}  \FloatBarrier

\vspace{1in}

\begin{table}[h]
	\caption{Group generators of 3D Laue groups \cite{ebbing2010design}.} \label{Table Laue generators}
	\centering
	\begin{tabular}{|c|c|}
	\hline
	Laue group & Group generators \\
	\hline
	$\mathcal C_i$     & $\bar{\mathbf I}$ \\
	$\mathcal{C}_{2h}$ & $\mathbf{C}_2, \bar{\mathbf I}$ \\ 
	$\mathcal{D}_{2h}$ & $\mathbf{C}_{2}, \mathbf C_{2x}, \bar{\mathbf I}$ \\
	$\mathcal{C}_{4h}$ & $\mathbf C_4, \bar{\mathbf I}$ \\
	$\mathcal{D}_{4h}$ & $\mathbf C_4,\mathbf C_{2x}, \bar{\mathbf I} $ \\
	$\mathcal{C}_{3i}  $ & $\mathbf C_3, \bar{\mathbf I}$ \\
	$\mathcal{D}_{3d}  $ & $\mathbf C_3, \mathbf C_{2x}, \ \bar{\mathbf I}$ \\	
	
	$ \mathcal{C}_{6h}   $ & $\mathbf C_6, \ \bar{\mathbf I} $\\
	$ \mathcal{D}_{6h}   $ & $\mathbf C_6,\ \mathbf C_{2x}, \ 
		\bar{\mathbf I} $\\
	$ \mathcal{T}_{h}$  & $\mathbf C_{2x},\ \mathbf C_{2y}, \mathbf Q^{2\pi/3}_p, \bar{\mathbf I}$ \\
	$ \mathcal{O}_{h}$ & $\mathbf C_{4x},\ \mathbf C_{2y}, \mathbf Q^{2\pi/3}_p, \bar{\mathbf I}$ \\ \hline
	\end{tabular}
\end{table}
	
Useful matrices for Table \ref{Table Laue generators}:
\begin{align*}
	\begin{array}{c}
			\mathbf I =\begin{bmatrix}
				1 & 0 & 0 \\
				0 & 1 & 0 \\
				0 & 0 & 1 
			\end{bmatrix},
			\bar{\mathbf I} =\begin{bmatrix}
			-1 & 0 & 0 \\
			0 & -1 & 0 \\
			0 & 0 & -1 
		\end{bmatrix},
		\mathbf{C}_{2}=\mathbf Q_{X_3}^\pi =\begin{bmatrix}
			-1 & 0 & 0 \\
			0 & -1 & 0 \\
			0 & 0 & 1 
		\end{bmatrix},
		\mathbf{C}_{2x}=\mathbf  Q_{X_1}^\pi =\begin{bmatrix}
			1 & 0 & 0 \\
			0 & -1 & 0 \\
			0 & 0 & -1 
		\end{bmatrix},\\
		\mathbf{C}_4=\mathbf Q_{X_3}^{\pi/2} = \begin{bmatrix}
			0 & 1 & 0 \\
			-1 & 0 & 0\\
			0 & 0 & 1
		\end{bmatrix},
		\mathbf{C}_3=\mathbf Q_{X_3}^{2\pi/3} = \begin{bmatrix}
			-1/2 & \sqrt 3 /2 & 0 \\
			-\sqrt 3 /2 & -1/2 & 0\\
			0 & 0 & 1
		\end{bmatrix},
		\mathbf{C}_6=\mathbf Q_{X_3}^{\pi/3} = \begin{bmatrix}
			1/2 & \sqrt 3 /2 & 0 \\
			-\sqrt 3 /2 & 1/2 & 0\\
			0 & 0 & 1
		\end{bmatrix}, \\
		\mathbf C_{2y}=\mathbf{Q}_{X_2}^\pi = \begin{bmatrix}
			-1&0&0\\
			0&1&0\\
			0&0&-1
		\end{bmatrix},
		\mathbf Q^{2\pi/3}_p = \begin{bmatrix}
			0 & 0&1\\
			1 & 0 & 0\\
			0 & 1 & 0
		\end{bmatrix},
		\mathbf C_{4x} =\mathbf{Q}_{X_1}^{\pi/2}= \begin{bmatrix}
			1 & 0 & 0\\
			0 & 0 & 1\\
			0 & -1 & 0
		\end{bmatrix}
	\end{array}
\end{align*}

\end{justify}
\end{document}